\documentclass[11pt]{amsart}

\usepackage{amscd,latexsym,amsthm,amsfonts,amssymb,amsmath,amsxtra,eucal}
\usepackage[all]{xy}
\usepackage[pdftex,bookmarks,colorlinks,pdfmenubar]{hyperref}

\usepackage{verbatim}
\usepackage{mathabx}
\usepackage{mathrsfs}
\usepackage{bbm}

\renewcommand\theequation{\thesection.\arabic{equation}}

\newcommand{\BC}{{\mathbb {C}}}

\newcommand{\BZ}{{\mathbb {Z}}}

\newcommand{\CP}{{\mathcal {P}}}

\newcommand{\Qlb}{\mathbb{\bar Q}_\ell}

\newcommand{\Ad}{{\mathrm{Ad}}}

\newcommand{\GL}{{\mathrm{GL}}}

\newcommand{\Hom}{{\mathrm{Hom}}}

\newcommand{\Spec}{{\mathrm{Spec}}}

\newcommand{\wt}{\widetilde}

\newcommand{\bs}{\backslash}

\def\diag{{\rm diag}}

\newtheorem{thm}{Theorem}[section]
\newtheorem{cor}[thm]{Corollary}
\newtheorem{lem}[thm]{Lemma}
\newtheorem{prop}[thm]{Proposition}

\newtheorem {ques/conj}[thm]{Question/Conjecture}

\newtheorem{defn}[thm]{Definition}

{\theoremstyle{definition}
\newtheorem{rmk}[thm]{Remark}}

\usepackage[numbers]{natbib}

\newcommand{\select}[1]{{\it{#1}}}

\begin{document}
\renewcommand{\theequation}{\arabic{equation}}
\numberwithin{equation}{section}

\title[Periods Of Deligne-Lusztig Characters Associated To Spherical Varieties]{Periods Of Deligne-Lusztig Characters Associated To Spherical Varieties}

\date{\today}

\author[Fang Shi]{Fang Shi}

\address{School of Mathematical Sciences, Xiamen University, Xiamen 361005, Fujian, P.R. China}

\email{shif@xmu.edu.cn}

\subjclass[2010]{Primary  20C33}

\maketitle

\begin{abstract}
In this paper, we calculate the periods of Deligne-Lusztig characters $R_{T,\chi}^G$ associated to $X^F$, where $X=G/H$ for a connected spherical subgroup $H$ of $G$. As an application, we deduce a formula that partially generalizes Lusztig's work on symmetric spaces and extends Reeder's formula in the complexity $0$ case. 
\end{abstract}

\section{Introduction}

\subsection{An overview}\label{overview}
We fix a finite field $\mathbb F_q$, where $q$ is a power of a prime $p$. We fix a prime $\ell$ that is different from $p$. We fix an identification $\Qlb \simeq \mathbb C$. We denote an algebraic closure of $\mathbb F_q$ by $\mathrm k$. 

For a scheme $X_0$ over $\mathbb F_q$, we will frequently denote its pullback to $\mathrm k$ by $X$. Then the geometric Frobenius $F$ acts naturally on $X$. We have a canonical bijection $X^F \simeq X_0(\mathbb F_q)$.

Let $G_0$ be a connected reductive group over $\mathbb F_q$ and let $H_0$ be a connected spherical subgroup of $G_0$. As mentioned above, we denote the pullback of $G_0$ (resp. $H_0$) to $\mathrm k$ by $G$ (resp. $H$).

Let $R_{T,\chi}^G$ be the virtual character of $G^F$ defined in \cite{DL} for a maximal $F$-stable torus $T$ and a character $\chi: T^F \to \Qlb^\times$. 
The main interest of this paper is to calculate 
$$\langle R_{T,\chi}^G, 1_{H^F}\rangle_{H^F}:=\frac{1}{|H^F|}\sum_{h\in H^F}R_{T,\chi}^G(h),$$
where we denote the trivial representation of $H^F$  by $1_{H^F}$.

Special cases, such as the symmetric spaces \cite{Lusym} and specific reductive subgroups \cite{R}, have
been studied extensively. While seeking a finite-field counterpart of \cite{WZ}, We express the number $\langle R_{T,\chi}^G, 1_{H^F}\rangle_{H^F}$ in terms of the geometric information of $G$ and $H$ in Theorem \ref{main} (under the sole assumption that the spherical subgroup $H$ is  connected). 

However, Theorem \ref{main}
 is a bit complicated to invoke here. In what follows, we introduce a refined version of Theorem \ref{main} under some assumption on $H$.

The assumption $\star$:
\begin{itemize}
    \item[] For any (not necessarily $F$-stable) Borel subgroup $B$ of $G$, the set $B(\mathrm k)\cap H(\mathrm k)$ equipped with the Zariski topology has an open dense subset $V$ consisting of semisimple elements.
\end{itemize}

We remark here that an element $g\in B(\mathrm k)\cap H(\mathrm k)$ is semisimple as an element in $B\cap H$ if and only if $g$ is a semisimple element in $G$ (\select{c.f.}, Theorem 9.18 and 9.19 of \cite{Mi}). And the words ``semisimple elements" in $\star$ should cause no confusion.

To elaborate our refined version of Theorem \ref{main} under the assumption $\star$, we introduce some notation.

Let $T$ and $\chi$ be as above. Let $B_T$ be a (not necessarily $F$-stable) Borel subgroup of $G$ containing $T$. Note that the group $B_T$ is not uniquely determined by $T$, but we eliminate this ambiguity in Definition \ref{defomeg}. Let $\mathrm d_T:B_T\to T$ be the obvious map witnessing $T$ as the reductive quotient of $B_T$ and providing a section for the inclusion $T\hookrightarrow B_T$. Let $\mathcal X:=G/H$.

\begin{defn}\label{d0-omega} (See Definition \ref{defomeg} for details.)
    Let $\Omega_T$ be the set of pairs $(c,X)$, where $c$ is a subscheme of $T$ and $X$ is a subscheme of $\mathcal X$ satisfying the following conditions:
    \begin{itemize}
        \item [i)] There exists $g\in G(\mathrm k)$ such that $c$ is a connected component of $\mathrm d_T(B_T\cap gHg^{-1})$, where we endow $\mathrm d_T(B_T\cap g H g^{-1})$ with the reduced scheme structure;
        \item[ii)] Let $\mathcal X^c$ be the subscheme of $\mathcal X$ fixed by $c$. Then $X$ is an irreducible component of $\mathcal X^c$ with $\dim X=\dim \mathcal X^c$.
    \end{itemize}
\end{defn}
Since the set $B_T(\mathrm k)\backslash G(\mathrm k)/H(\mathrm k)$ is  finite, there are only finitely many $c$ as in i) of Definition \ref{d0-omega}. We see from Definition \ref{defomeg} that $\Omega_T$ does not depend on the choice of $B_T$, justifying the notation. And $\Omega_T$ is a finite set endowed with a natural action of the Frobenius endomorphism given by $(c,X)\mapsto(F(c),F(X))$. 

For a connected algebraic  group $K$ over $\mathbb F_q$, we denote the $\mathbb F_q$-rank  of $K$ by $\sigma(K)$. For a subscheme $Y$ of $G$, we denote the centralizer of $Y$  in $G$ by $C_G(Y)$.
For an algebraic group $K$, we denote its identity component by $K^\circ$.
\begin{defn}\label{d0-tphi}
    (See Definition \ref{deftphi} for details.) Let $\omega=(c,X)\in \Omega_T^F$, \select{i.e.}, the subscheme $c$ ($X$, resp.) of $T$ ($\mathcal X$, resp.) is $F$-stable. We set $$
    \sigma_\omega:=\sigma(C_G(c)^\circ) \quad\text{and}\quad
    t_{\omega,\chi}:=\frac{1}{|c^F|}\sum_{s\in c^F}\chi(s).
    $$
\end{defn}
Since the subscheme $c$ of $T$ in Definition \ref{d0-tphi} is an $F$-stable component of an algebraic subgroup of $T$, it is conspicuous that $t_{\omega,\chi}$ is either zero or a root of unity.
The following is our refined main theorem under the assumption $\star$. (See Theorem \ref{emmain} for details.)

\begin{thm}
    Suppose that the assumption $\star$ is fulfilled. Let $T$  and $\chi$ be as above. We have 
    $$
    \langle R_{T,\chi}^G,1_{H^F}\rangle_{H^F}=\sum_{\omega\in \Omega_T^F}(-1)^{\sigma(T)+\sigma_\omega}t_{\omega,\chi}.
    $$
\end{thm}

\begin{rmk}\label{rm-example}
    The assumption $\star$ seems to be artificial. However, We will see that a vast number of interesting examples satisfy this assumption. Here are some of them.
    \begin{itemize}
        \item [(i)] $H$ is an $F$-stable parabolic subgroup of $G$; (See Section \ref{e-para}.)
        \item [(ii)] $H$ is the connected subgroup of $G$ fixed by an involution; (See Section \ref{e-sym}.)
        \item [(iii)] The diagonal embedding $H=G_1\hookrightarrow G_1\times G_2=G$, where $G_1$ is a complexity $0$ reductive subgroup of the reductive group  $G_2$ in the sense of \cite{R}; (See Section \ref{e-Reeder}, where we include finite-field analogs of the basic cases of Bessel models in the sense of \cite{GGP}.)
        \item [(iv)] Some finite-field analogs of models discussed in \cite{WZ}, \select{e.g.}, the pair $(\mathrm U_4\times \mathrm U_2,\mathrm U_2\times \mathrm U_2)$. (See Section \ref{exmain}.)
    \end{itemize}
\end{rmk}

In particular, Theorem \ref{emmain} agrees with Theorem 3.3 of \cite{Lusym} when $H$ is the connected subgroup fixed by an involution. Also, Theorem \ref{emmain} reproduces results concerning the complexity $0$ case of \cite{R} (\select{e.g.}, the restriction problem for $\mathrm {SO}_{2n}\hookrightarrow \mathrm{SO}_{2n+1}$). 
 
\subsection{Structure of this paper}Section \ref{s-transpsemi} and Section \ref{s-multiplici} are dedicated to deduce Theorem \ref{main}.
Here we outline the strategy to calculate the relevant branching number. 
Inspired by \cite{R}, we replace the Frobenius $F$ with powers, giving rise to the  function 
$\mathrm P(\nu,T,\chi)$ for $\nu$ varying over an arithmetic progression $\mathcal P$. Moreover, we have $\mathrm P(1,T,\chi)= \langle R_{T,\chi}^G, 1_{H^F}\rangle_{H^F}$. (See Definition \ref{dfp} and Definition \ref{dtp} for  precise definitions.) Then we show $\mathrm P(\nu,T,\chi)$ has a finite limit  as $\nu \to \infty$. As a result of  Lemma \ref{const}, the function $\mathrm P(\nu, T,\chi)$  is constant as $\nu$ varies over $\mathcal P$. 
 
Similar approaches are taken in  \cite{LMS}.

 For a positive integer $d$, we set $\mathcal P_{d}:=\{1+\nu d\}_{\nu \in \mathbb N}$ to be the arithmetic progression starting from $1$ with the gap $d$.

  The following is  a corollary of Theorem \ref{main}. 

\begin{thm}
For any $F$-stable maximal torus $T$ of $G$, there exists a positive integer $d_T$ such that for any character $\chi: T^F \to \Qlb^\times$, the function $\mathrm P(\nu,T,\chi)$
is a constant as $\nu$ varies over $\mathcal P_{d_{T}}$.
\end{thm}

In Section \ref{s-semisimpleele}, we simplify the complex multi-indices introduced in Theorem \ref{main} under the assumption $\star$. The main result of Section \ref{s-semisimpleele} is Theorem \ref{emmain}.

In Section \ref{s-examp}, we consider certain examples to which Theorem \ref{emmain} applies. (However, we use Theorem \ref{main} to deal with $(\mathrm U_4\times \mathrm U_2,\mathrm U_2\times \mathrm U_2)$. See Remark \ref{r-swmain}.) We incorporate some results from \cite{Lusym} (when the relevant subgroup  is connected) and \cite{R}
(when the complexity is $0$) into the framework of Theorem \ref{emmain}, as promised in the previous subsection.

In the appendix, we recall some results of \cite{L5}. (The corresponding results are only used in Definition \ref{dtp} and Remark \ref{spc}.)
\subsection{Conventions}\label{conv}

For a scheme $X$ of finite type over $\mathrm k$, we denote its Krull dimension by $\dim X$. Set $\dim X=-\infty$ if $X$ is the empty scheme. We denote the derived category of bounded complexes of constructible $\Qlb$-sheaves on $X$ by $D(X,\Qlb)$.


 For a finite group $K$, we will frequently denote $\dim \Hom_K(\tau,\pi)$ by $\langle \tau,\pi\rangle_K$ for a pair of $\Qlb$-representations $(\tau,\pi)$ of $K$. We will abuse the notation by denoting the linear extension of $\langle -,-\rangle_K$ to virtual characters of $K$ again by $\langle-,-\rangle_K$. We denote the trivial representation of $K$ by $1_K$.
 
 Let $\mathbb Z_+$ be the set of positive integers.

For a positive integer $d$, we set $\mathcal P_{d}:=\{1+\nu d\}_{\nu \in \mathbb N}$ to be the arithmetic progression starting from $1$ with the gap $d$.

For a semisimple element $s\in G(\mathrm k)$, we denote the identity component of the centralizer $C_{G}(s)$ by $G_s$.


We often identify a scheme $X_0$ over $\mathbb F_q$ with its pullback $X$ to $\mathrm k$ equipped with the Frobenius endomorphism $F:X \to X$. This should cause no confusion.

Assume that $q$ is large enough for  Dynkin diagrams of all algebraic groups of the form  $G_s$  in the sense of Theorem 1.14 of \cite{L5} where $s\in G(\mathrm k)$ is semisimple.  

For an endomorphism $E$ of a $\Qlb$-linear space $V$ of finite dimension, we denote the trace of $E$ by $\mathrm{Tr}(E,V)$. 

For a scheme $X$ over $\mathrm {k}$ and $x \in X(\mathrm {k})$, we sometimes write $x \in X$ to ease the notation. 

For an arithmetic progression $\mathcal P$, we write $\mathcal P \ni \nu \to \infty$ to indicate that the value of the variable $\nu$ in $\mathcal P$ tends to infinity.

Let $X_0$ and $X$ be as in the first paragraph of Section \ref{overview}. For a Weil sheaf $(\mathscr F,a:F^* \mathscr F \widetilde{\to} \mathscr F)$ on  $X$, we denote the trace of the endomorphism of $R\Gamma_c(X,\mathscr F)$ induced by $\mathscr F \stackrel{adj}{\to }F_*F^* \mathscr F \stackrel{F_*(a)}{\to} F_* \mathscr F$ by $\mathrm{Tr}(F,R\Gamma_c(X,\mathscr F))$. Similarly we define the trace $\mathrm{Tr}(F,\mathscr F_x)$ for an $F$-stable point $x\in X^F$. 
\subsection*{Acknowledgment}
Some portions of this paper constitute a significant  part of the author's thesis.
The author is indebted to Dongwen Liu, who suggested the author consider the
finite-field analogs of several models  in \cite{WZ}. Also, the author is grateful for valuable discussion with Dongwen Liu and Zhicheng Wang.

\section{Transporters of semisimple elements}\label{s-transpsemi}
We fix an $F$-stable maximal torus $S$ of $H$ and an $F$-stable maximal torus $T$ of $G$ in this section.  
\subsection{Transporters with reduced scheme structure}
For any $r\in G(\mathrm k)$, let $\mathrm T_\circ(r,H)$ be the closed subscheme of $G$ representing the  functor via Yoneda embedding sending an algebra $R$ over $\mathrm k$ to
$$
\{g\in G(R) :  g^{-1} r g \in H(R) \},
$$
where we view $r$ as an $R$-point via the inclusion $G(\mathrm k) \hookrightarrow G(R)$.
We have a pullback square
$$
\xymatrix{
\mathrm T_\circ(r,H) \ar[r]\ar[d]  & H\ar[d]\\
G\ar[r]&G\\
}
$$
where the right vertical arrow is the inclusion and the bottom arrow sends $g$ to $g^{-1} r g$. Let $\mathrm T(r,H)$ be the closed subscheme of $\mathrm T_\circ(r,H)$ with the same underlying topological space and the reduced scheme structure.

\subsection{$\mathrm k$-points of $\mathrm T(t,H)$}
 For $t\in T(\mathrm k)$, we set $N(t,S,T):=\{g\in G(\mathrm k):t \in gSg^{-1}(\mathrm k), gSg^{-1}\subset T\}$.
\begin{lem}\label{kpt} Let $t\in T(\mathrm k)$.
The set of the $\mathrm k$-points of $\mathrm T(t,H)$ is the image of $G_t(\mathrm k)\times N(t,S,T)\times H(\mathrm k)$ under the multiplication of $G(\mathrm k)$.
\end{lem}
\begin{proof}(See Section \ref{conv} for the notation $G_t$.)
It is clear that the image of the above multiplication is a subset of $\mathrm T(t,H)(\mathrm k)$.
Suppose that we have $g^{-1}tg \in H(\mathrm k)$. Since $t$ is semisimple and $S$ is a maximal torus of $H$, we can take some $h\in H(\mathrm k)$ so that $h^{-1}g^{-1}tgh \in S(\mathrm k)$. Hence $ghSh^{-1}g^{-1}\subset G_t$. Since $T$ is a maximal torus of $G_t$, we can take $\wt g \in G_t(\mathrm k)$ so that 
$\wt g g h Sh^{-1}g^{-1}\wt g^{-1} \subset T$. And it is easy to see $\wt g g h \in N(t,S,T)$, as desired.
\end{proof}

\subsection{A partition of $T$}
Let $\mathcal J(T,S)$ be the set of subtori of $T$ that are $G(\mathrm k)$-conjugates of $S$. 
Here we see that different choices of $S$ give the same $\mathcal J(T,S)$. Hence we set $\mathcal J(T):=\mathcal J(T,S)$ for some $F$-stable maximal torus $S$ of $H$.
For any nonempty subset $\jmath \subset J(T)$, let $T_\jmath$ be the following locally closed subscheme of $T$ 
$$\bigcap\limits_{K\in\jmath}K-\bigcup_{R\in \mathcal J(T) -\jmath}R$$
equipped with the reduced scheme structure. If $\jmath$ is the empty set, we set $T_\jmath=T-\bigcup\limits_{K \in \mathcal J(T)} K$. %
Let $I(T)$ be an index set for the set of  subgroups $\{G_t: t \in T(\mathrm k)\}$. For $\iota \in I(T)$, let $G_\iota$ be the corresponding connected centralizer, and we set $T_\iota$ to be the locally closed reduced  subscheme of $T$ whose set of $\mathrm k$-points is
$$\{t\in T(\mathrm k):G_t =G_\iota\}.
$$
We set $T_{\jmath,\iota}=T_\jmath \cap T_\iota$ with the reduced scheme structure for $\jmath \subset \mathcal J(T)$ and $\iota \in I(T)$. We denote the power set of $\mathcal J(T)$ by $2^{\mathcal J(T)}$.

\begin{rmk}\label{rmtos}
   For $(\jmath,\iota)\in 2^{\mathcal J(T)}\times I(T)$,  we see from the definition of $T_{\jmath,\iota}$ that the scheme $T_{\jmath,\iota}$ is an open subscheme of an algebraic group. (Indeed, it is an open subscheme of the algebraic group $\dot{T}_{\jmath,\iota}$ introduced in Definition \ref{tcirc}.)
\end{rmk}

\begin{prop}\label{st}
The set $\{T_{\jmath,\iota}  \}_{(\jmath,\iota)}$ of reduced locally closed subschemes of $T_\iota$ indexed by the set $2^{\mathcal J(T)} \times I(T)$ forms a partition of $T$. For each pair $(\jmath,\iota)\in 2^{\mathcal J(T)}\times I(T)$, the following statement holds:\\
For any $t_1,t_2 \in T_{\jmath,\iota}$, we have $$\mathrm T(t_1,H)=\mathrm T(t_2,H).$$
\end{prop}
\begin{proof}
The first statement is clear. We prove the last statement.
Since $\mathrm T(t_1,H)$ and $\mathrm T(t_2,H)$ are reduced subschemes of $G$, we need to show
$$\mathrm T(t_1,H) (\mathrm k)=\mathrm T(t_2,H)(\mathrm k).
$$
Thanks to Lemma \ref{kpt}, it remains to show $N(t_1,S,T)=N(t_2,S,T)$. We see $N(t_1,S,T)=\bigcup\limits_{K\in \jmath} N_K=N(t_2,S,T)$, where $N_K:=\{g\in G(\mathrm k):g Sg^{-1}=K\}$. This completes the proof.
\end{proof}

\begin{rmk}
The Frobenius $F$ acts naturally on $\mathcal J(T)$, which gives rise to an action on $2^{\mathcal J(T)}$. This action satisfies that $F(T_\jmath)=T_{F(\jmath)}$ for $\jmath \in 2^{\mathcal J(T)}$. Similarly, the Frobenius $F$ acts naturally on $I(T)$.  We have $F(T_\iota)=T_{F(\iota)}$ and $F(T_{\jmath,\iota})=T_{F(\jmath),F(\iota)}$ for $\iota \in I(T)$ and $\jmath \in \mathcal J$. The set of $F$-invariant elements of $2^{\mathcal J(T)}$ (resp. $I(T)$) is denoted by $2^{\mathcal J(T),F}$ (resp. $I(T)^F$).
\end{rmk}

\begin{rmk}\label{r-gpartition}
    The assumption that $T$ and $S$ are $F$-stable guarantees that the sets $\mathcal J(T)$ and $I(T)$ possess natural actions of the Frobenius $F$.
    Suppose that $T'$ is a (not necessarily $F$-stable) maximal torus of $G$ and $S'$ is a (not necessarily $F$-stable) maximal torus of $H$. The sets $\mathcal J(T')$ and $I(T')$ can be likewise defined. (However, we do not pursue an action of the Frobenius on $\mathcal J(T')$ and $I(T')$.)
\end{rmk}

\section{The Multiplicity Formula}\label{s-multiplici}
\subsection{Functions of geometric type} We introduce the following notion, which is used in \cite{LMS}.
\begin{defn} \label{gt}
Let  $\CP\subset \BZ_+$ be an arithmetic progression. A function $M: \CP \to \BC $ is said to be of geometric type if it is of the form
\[
M(\nu) = \frac{\sum^k_{i=1} a_i \alpha_i^\nu}{\sum^l_{j=1} b_j \beta_j^\nu } ,\quad \nu\in\CP,
\]
where $a_i, \alpha_i, b_j, \beta_j \in\BC$, and the denominator is nonzero for every $\nu\in \CP$. If the denominator is a nonzero constant, we say $M$ is of trace type.  
\end{defn}
We have the following elementary lemma, see Lemma 2.2 of \cite{LMS}. 

\begin{lem} \label{const}
Let $M$ be a function of geometric type defined on an arithmetic progression $\CP\subset \mathbb Z_+$. If $M$ is integer-valued and  has a finite limit $L\in \mathbb C$ as $\nu\to\infty$ along $\CP$, then 
$M$ is a constant function taking the value $L$. 
\end{lem}

\subsection{Summation on $G^F$} In this subsection, we reformulate Section 5.1 of \cite{R}. 

Suppose we have a function $f: G^F \to \mathbb C$, invariant  under conjugation by $G^F$, with the property that
\begin{itemize}\label{fp}
\item if $g \in G^F$ has Jordan decomposition $g=su$, then $f(g)=0$ unless the conjugacy class $\Ad(G^F)\cdot s$ meets $T^F$. (Here, the element $s\in G$ is semisimple and $u\in G$ is unipotent.)
\end{itemize}
\begin{defn}\label{d-barnugs}
Recall that for $s\in T(\mathrm k)$, we denote the identity component of the centralizer $C_G(s)$ by $G_s$. For $r\in G(\mathrm k)$, we set $N(r,T)$ to be the reduced subscheme of $G$ whose set of  $\mathrm k$-points is $\{g\in G(\mathrm k):g^{-1}r g \in T(\mathrm k)\}$. 
We set $\bar N(s,T):=G_s\backslash N(s,T)$.
 Let $\mathcal U(G_s^F)$ be the set of $\Ad(G_s^F)$-orbits of unipotent elements in $G^F_s$. 
 \end{defn}
 For $s\in T^F$, we see that the fixed point set  of $\bar N(s,T)$ under the Frobenius $F$ is $G_s^F \bs N(s,T)^F$. The following is (5.2) of \cite{R}.
\begin{prop}\label{sf}Let $f$ be as above.
The following equation holds:
 $$
\frac{1}{|G^F|}\sum_{g\in G^F}f(g)=\sum_{s\in T^F}\frac{1}{|\bar N(s,T)^F|}\sum_{[u]\in \mathcal U(G_s^F)} \frac{1}{|C_{G_s}(u)^F|}f(su).
$$
\end{prop}
Recall that we have a partition of $T$ indexed by $2^{\mathcal J(T)}\times I(T)$, which gives rise to a partition of $T^F$ indexed by $2^{\mathcal J(T),F}\times I(T)^F$. Note that for $\iota \in I(T)$ and $s_1,s_2\in T_\iota(\mathrm k)$, we have $G_{s_1}=G_{s_2}$ and $N(s_1,T)=N(s_2,T)$. For $\iota\in I(T)$, we denote $G_\iota:=G_s$ (resp. $N(\iota,T):=N(s,T)$ and $\bar N(\iota,T):=\bar N(s,T)$) for some $s\in T_\iota(\mathrm k)$. 

For $\iota \in I(T)^F$, we verify that $G_\iota$ and $N(\iota,T)$ are $F$-stable. And we denote the fixed point set  of $\bar N(\iota,T)$ under the Frobenius $F$ by $\bar N(\iota,T)^F$.

 We have the following corollary.
\begin{cor}\label{fst}Let $f$ be as above. Then we have the following:
$$
\frac{1}{|G^F|}\sum_{g\in G^F}f(g)=\sum_{\substack{\jmath \in 2^{\mathcal J(T),F} \\ \iota \in I(T)^F}}\sum_{s\in T_{\jmath,\iota}^F}\frac{1}{|\bar N(\iota,T)^F|}\sum_{[u]\in \mathcal U(G_\iota^F)} \frac{1}{|C_{G_\iota}(u)^F|}f(su).$$
\end{cor}

\subsection{Deligne-Lusztig characters} \label{dl}
Recall that $T$ is an $F$-stable maximal torus of $G$. Let 
$\chi: T^F \to \Qlb^\times$ be a character. The virtual character $R_{T,\chi}^G$ is defined in \cite{DL}. Let $s$ be a semisimple element of $G^F$. For a unipotent element $u\in G_s^F$, we have the following equation (this is a reformulation of Theorem 4.2 of \cite{DL}, see also (4.1) of \cite{R}):

\begin{equation}\label{fdl}
R_{T,\chi}^G(su)= \sum_{\bar \gamma \in \bar N(s,T)^F} \chi(\gamma^{-1}s \gamma) Q_{\gamma T\gamma^{-1}}^{G_s}(u).
\end{equation}
Here $ Q_{\gamma T\gamma^{-1}}^{G_s}:=R_{\gamma T\gamma^{-1},1}^{G_s}$ is the Green function.
\begin{rmk}\label{fpdl}
We see that 
the function $f_0 \cdot R_{T,\chi}^G$ has the property displayed in Section \ref{fp} for any $F$-stable maximal torus $T$ of $G$, any character $\chi: T^F \to \Qlb^\times$ and any virtual character $f_0$ of $G^F$.
\end{rmk}

For future use, we define $R_{T,\chi}^{G,\nu}$ to be the virtual character $R_{T,\chi\circ \mathrm N^\nu}^G$ of $G^{F^\nu}\simeq G_0(\mathbb F_{q^\nu})$ for $\nu \in \mathbb Z_+$. Here we view $T$ as an $F^\nu$-stable maximal torus of $G\simeq \left(G_0 \times_{\Spec (\mathbb F_q)} \Spec(\mathbb F_{q^\nu})\right)\times_{\Spec (\mathbb F_{q^\nu})} \Spec(\mathrm k)$ and $\mathrm N^\nu: T^{F^\nu} \to T^F$ is the norm map. Similarly we write  $Q_{\gamma T\gamma^{-1}}^{G_s,\nu}$ for $R_{\gamma T\gamma^{-1},1}^{G_s,\nu}$.

\subsection{Flag variety}\label{fva}
Fix an $F$-stable Borel subgroup $B$ of $G$. Then the flag variety $\mathcal B$ that classifies Borel subgroups of $G$ is isomorphic to $G/B$. Throughout this section, we fix this identification.
For $g\in G$, let $\mathcal B_g:=\{xB \in G/B:x^{-1}gx \in B\}$ be the subscheme of the flag variety $\mathcal B$ of $G$. 
 Suppose that $g=su$ is the Jordan decomposition. Then we have $\mathcal B_g= \mathcal B_s \cap \mathcal B_u$ as subschemes of $\mathcal B$. Fix $\iota\in 2^{\mathcal J(T)}\times I(T)$. We verify that for $s_1,s_2\in T_\iota(\mathrm k)$ we have $\mathcal B_{s_1}=\mathcal B_{s_2}$. Hence for unipotent $u \in G_\iota$, we have $\mathcal B_{s_1u}=\mathcal B_{s_2u}$, which we denote by $\mathcal B_{\iota,u}$. Let $d_{\iota,u}:=\dim \mathcal B_{\iota,u}$.

\subsection{Periods of Deligne-Lusztig characters}
 Recall that for $(\jmath,\iota)\in 2^{\mathcal J(T)}\times I(T)$ and $s_1,s_2 \in T_{\jmath,\iota}(\mathrm k)$, we have $\mathrm T(s_1,H)=\mathrm T(s_2,H)$ by Proposition \ref{st}. For  $(\jmath,\iota)\in 2^{\mathcal J(T)}\times I(T)$ we set $\mathrm T(\jmath,\iota):=\mathrm T(s,H)$ for some $s\in T_{\jmath,\iota}(\mathrm k)$. Note that for $(\jmath,\iota) \in 2^{\mathcal J(T),F}\times I(T)^F$, the reduced scheme $\mathrm T(\jmath,\iota)$ is $F$-stable since its set of  $\mathrm k$-points is $F$-stable by Lemma \ref{kpt}. For $(\jmath,\iota) \in 2^{\mathcal J(T),F}\times I(T)^F$ and unipotent $u \in G_\iota ^F$, we set $\mathrm T(\jmath,\iota,u):= \mathrm T(\jmath,\iota)\cap \mathrm T(u,H)$ to be  the ($F$-stable) reduced subscheme of $G$ whose set of $\mathrm k$-points  is $\mathrm T(\jmath,\iota)(\mathrm k) \cap \mathrm T(u,H)(\mathrm k)  $.

We introduce the following symbols:
\begin{defn}\label{d-mathcalC}
Let $X_0$ be a scheme of finite type over $\mathbb F_q$ of dimension $d$. Let $X$  be the pullback of $X_0$ to $\mathrm k$. We define $\mathcal C(X)$ to be the number of $d$-dimensional $F$-stable irreducible components of $X$.  
\end{defn}

\begin{defn} \label{dfp}
Let $T$ be a $F$-stable maximal torus of $G$. Let $\chi: T^F \to \Qlb$ be a character.
Recall the virtual character $R_{T,\chi}^{G,\nu}$ of $G^{F^\nu}$ introduced in the end of Section \ref{dl}.
For $\nu \in \mathbb Z_+$,
we set $$\mathrm P (\nu,T,\chi):=\langle R_{T,\chi}^{G,\nu},1_{H^{F^\nu}}\rangle_{H^{F^\nu}}.$$
\end{defn}

\begin{defn}\label{tcirc}
 For $\iota \in I(T)$ we set ${\dot{T}}_\iota$ to be the reduced subgroup of $T$ fixed by $G_\iota$.
 For $\jmath \in 2^{\mathcal J(T)}$, we set ${\dot{T}}_\jmath$ to be the reduced subgroup of $T$ whose set $\mathrm k$-points is $\bigcap\limits_{K\in \jmath} K(\mathrm k)$.
 For $\iota \in I(T)$ and $\jmath \in 2^{\mathcal J(T)}$, we set ${\dot{T}}_{\jmath,\iota}$ to be the reduced subgroup of $T$ whose set of $\mathrm k$-points is ${\dot{T}}_{\iota}(\mathrm k) \cap {\dot{T}}_\jmath(\mathrm k)$. For a 
 subset $\Theta$ of $I(T)$, we set ${\dot{T}}_{\jmath,\Theta}$ to be $\bigcap\limits_{\iota \in \Theta} {\dot{T}}_{\jmath,\iota}$ to be the closed subgroup of $T$ equipped with the reduced subscheme structure.

\end{defn}

\begin{defn}\label{dtp}
For each $F$-stable maximal torus $T$ of $G$,
we fix an integer $d_T \in \mathbb Z_+$, so that the following conditions hold
:
\begin{itemize} 
\item $F^{d_T}$ acts trivially on $I(T)$ and $\mathcal J(T)$;
\item For $\iota \in I(T)^F$ and unipotent $u\in G_{\iota}^F$, set $A_\iota(u)$ to be the component group of $C_{G_\iota}(u)$. Then $c_u:=\frac{d_T}{|A_\iota(u)|}$ is an integer, and  the endomorphism $F^{c_u}$ acts trivially on $A_\iota (u)$;
\item For $\iota \in I(T)^F$, the endomorphism $F^{d_T}$ acts trivially on $\bar N(\iota, T)$;
\item For $\iota \in I(T)^F$ and $\jmath \in 2^{\mathcal J(T),F}$, the endomorphism $F^{d_T}$ acts trivially on the set of irreducible components of $\mathrm T(\jmath,\iota)$;
\item For $\iota \in I(T)^F$ and $\jmath \in 2^{\mathcal J(T),F}$, the endomorphism $F^{d_T}$ acts trivially on the set of irreducible components of $T_{\jmath,\iota}$;
\item For any unipotent element $u\in G^F$, the endomorphism $F^{d_T}$ acts trivially on the set of irreducible components of $\mathrm T(u,H)$;
\item For $(\jmath,\iota)\in 2^{\mathcal J(T),F}\times I(T)^F$ and  unipotent $u\in G_\iota^F$, the endomorphism $F^{d_T}$ acts trivially on the set of irreducible components of $\mathrm T(\jmath,\iota,u)$;
\item For any $F$-stable subset $\Theta$ of $I(T)$ and $\jmath \in 2^{\mathcal J(T),F}$, the endomorphism $F^{d_T}$ acts trivially on the set of irreducible components of ${\dot{T}}_{\jmath,\Theta}$;
\item For any $F$-stable subset $\Theta$ of $I(T)$ and $\jmath \in 2^{\mathcal J(T),F}$, the number $\frac{d_T}{\mathcal C({\dot{T}}_{\jmath,\Theta})}$ is an integer;
\item Take any $\iota \in I(T)^F$, $\gamma \in N(\iota,T)^F$ and unipotent $u\in G_\iota^F$.
There are  functions  $Q_{\iota,u,\gamma}$ of trace type in the sense of Definition \ref{gt}  satisfying that $Q_{\gamma T\gamma^{-1}}^{G_\iota,\nu}(u)=Q_{\iota,u,\gamma}(\nu)$ for $\nu \in \mathcal P_{d_T}$,  and the limit (see Subsection \ref{fva} for the definition of $d_{\iota,u}$) 
$$\lim \limits_{\mathcal P_{d_T} \ni \nu \to \infty} \frac{Q_{\iota,u,\gamma}(\nu)}{q^{\nu d_{\iota,u}}}$$
exists.  We denote the above limit by $q_{\iota,u,\gamma}$.
\end{itemize}
\end{defn}
\begin{rmk}\label{spc}
An integer $d$ that satisfies the last condition of Definition \ref{dtp} can be derived from Theorem \ref{dls} and Proposition \ref{lmq}. In particular, if $u=1$, then we have $q_{\iota,u,\gamma}=(-1)^{\sigma(G_\iota)-\sigma(T)}$ by Theorem 7.1 of \cite{DL}. 
\end{rmk}

\begin{rmk}
    Definition \ref{dtp} seems to be clumsy. All conditions listed are to ensure that the index sets of the (outer three) summations on the right-hand side of Equation (\ref{tff}) remain unchanged and the inner sum there has a finite limit as $\mathcal P_{d_T}\ni\nu \to \infty$ .
\end{rmk}

We will prove the following proposition in the next subsection (see also Section \ref{sadlc} for an alternative proof):
\begin{prop}\label{jgt}
Fix an $F$-stable torus $T$ of $G$ and a character $\chi :T^F \to \Qlb^\times$. 
Then the function $M(\nu):=\mathrm P(\nu,T,\chi)$ is of geometric type with respect to the arithmetic progression $\mathcal P_{d_T}$ in the sense of Definition \ref{gt}.
\end{prop}

By Frobenius reciprocity, we have the following:
\begin{equation}
\mathrm P(\nu,T,\chi)=\langle R_{T,\chi}^{G,\nu},\mathrm{Ind}_{H^{F^\nu}}^{G^{F^\nu}} 1_{H^{F^\nu}}\rangle_{G^{F^\nu}}=\frac{1}{|G^{F^\nu}|}
\sum_{g\in G^{F^\nu}} R_{T,\chi}^{G,\nu}(g) \cdot \frac{|\mathrm T(g,H)^{F^\nu}|}{|H^{F^\nu}|}.
\end{equation}
By Remark \ref{fpdl}, Equation (\ref{fdl})  and Corollary \ref{fst}, we see
\begin{equation}\label{tff}
\mathrm P(\nu,T,\chi)=\sum_{\substack{\jmath \in 2^{\mathcal J(T),F^\nu} \\ \iota \in I(T)^{F^\nu}}}\sum_{\bar \gamma \in \bar N(\iota,T)^{F^\nu}}\sum_{[u]\in \mathcal U(G_\iota^{F^\nu})}\sum_{s\in T_{\jmath,\iota}^{F^\nu}}\frac{|\mathrm T(su,H)^{F^\nu}|\cdot   Q_{\gamma T\gamma^{-1}}^{G_\iota,\nu}(u)}{|\bar N(\iota,T)^{F^\nu}|\cdot|C_{G_\iota}(u)^{F^\nu}|\cdot|H^{F^\nu}|}\cdot\chi\circ \mathrm N^ \nu(\gamma^{-1}s \gamma).
\end{equation}

Let $g=su \in G^F$ be the Jordan decomposition. Assume that $s\in T_{\jmath,\iota}$ for $(\jmath,\iota)\in 2^{\mathcal J(T),F}\times I(T)^F$. It is clear that $$|\mathrm T(g,H)^F|=|\mathrm T(s,H)^F\cap \mathrm T(u,H)^F|=|\left(\mathrm T(\jmath,\iota)\cap \mathrm T(u,H)\right)^F|=|\mathrm T(\jmath,\iota,u)^F|.$$

By our assumptions on $d_T$, we see that for $\nu \in \mathcal P_{d_T}$, the index sets for the outer three sum in Equation (\ref{tff}) remain unchanged. Hence we get the following:
\begin{lem} \label{meq}
Fix a $F$-stable torus $T$ of $G$ and a character $\chi :T^F \to \Qlb^\times$. For $\nu \in \mathcal P_{d_T}$, we have
$$\mathrm P(\nu,T,\chi)=\sum_{\substack{\jmath \in 2^{\mathcal J(T),F} \\ \iota \in I(T)^{F}}}\sum_{\bar \gamma \in \bar N(\iota,T)^{F}}\sum_{[u]\in \mathcal U(G_\iota^{F})}\frac{|\left(\mathrm T(\jmath,\iota)\cap \mathrm T(u,H)\right)^{F^\nu}|\cdot   Q_{\gamma T\gamma^{-1}}^{G_\iota,\nu}(u)}{|\bar N(\iota,T)^{F}|\cdot|C_{G_\iota}(u)^{F^\nu}|\cdot|H^{F^\nu}|}\sum_{s\in T_{\jmath,\iota}^{F^\nu}}\chi\circ \mathrm N^ \nu(\gamma^{-1}s \gamma).
$$
\end{lem}

\subsection{Periods are of geometric type for $\mathcal P_{d_T}$} In this subsection, we prove Proposition \ref{jgt}.
We need the following lemma, which is a trivial instance of the Grothendieck trace formula:
\begin{lem}\label{sfc}
Let $X_0$ be a separated scheme of finite type over $\mathbb F_q$. Let $X$ be the pull-back of $X_0$ to $\mathrm k$. Suppose that $\mathscr L \in D(X_0,\Qlb)$. Then the function
$$ M(\nu)=\sum_{x \in X^{F^\nu}} \mathrm{Tr}(F^\nu,\mathscr L_x)
$$
is of trace type for the arithmetic progression $\mathbb Z_+$. 
 In particular, the function $M$ is of trace type for any arithmetic progression $\mathcal P \subset \mathbb Z_+$.
\end{lem}

\begin{rmk}\label{cst}
Let $R_0$ be a torus over $\mathbb F_q$. Let $R$ be its pullback to $\mathrm k$.
The map $L^{R_0}:R_0\to R_0$ defined by $t \mapsto F(t)^{-1}t$ is finite \'etale. The local system $L_* ^{R_0}\Qlb$ is a direct sum
$\oplus_\eta \mathscr L_\eta$, where $\eta$ varies over the set of the  characters of $R^F$, and the local system $\mathscr L_\eta$ satisfies that $\mathrm{Tr}(F^\nu,(\mathscr L_\eta)_x)=\eta \circ \mathrm N^\nu(x)$ for any $\nu \in \mathbb Z_+$ and $x \in R^{F^\nu}$. ($\mathrm N^\nu: R^{F^\nu}\to R^F$  is the norm map.)
\end{rmk}

To prove Proposition \ref{jgt}, it suffices to show the following proposition.
\begin{prop}
The following functions are of trace type for $\mathcal P_{d_T}$ in the sense of Definition \ref{gt} (in the following, we define the functions for $\nu \in \mathcal P_{d_T}$):
\begin{itemize}
\item[(1)] the functions $M_1(\nu)=|\left(\mathrm T(\jmath,\iota)\cap T(u,H)\right)^{F^\nu}|=|\mathrm T(\jmath,\iota,u)|$ for $(\jmath,\iota) \in 2^{\mathcal J(T),F} \times I(T)^F$ and unipotent $u \in G_\iota^F$;
\item[(2)] the functions $M_2(\nu)= |C_{G_\iota}(u)^{F^\nu}|$ for $\iota \in I(T)^F$ and unipotent $u \in G_\iota^F$;
\item[(3)] the function $M_3(\nu)=|H^{F^\nu}|$;
\item[(4)] the functions $M_4(\nu)=\sum_{s\in T_{\jmath,\iota}^{F^\nu}}\chi\circ \mathrm N^ \nu(\gamma^{-1}s \gamma)$ for $(\jmath,\iota) \in 2^{\mathcal J(T),F} \times I(T)^F$ and $\gamma \in N(\iota,T)^F$;
\item[(5)] the functions $M_5(\nu)=Q^{G_\iota,\nu}_{\gamma T \gamma^{-1}}(u)$ for $\iota \in I(T)^F$ and unipotent $u \in G_\iota^F$.
\end{itemize}
\end{prop}
\begin{proof}
For cases (1),(2) and (3), we apply Lemma \ref{sfc} to the corresponding schemes over $\mathbb F_q$ and the constant sheaves of rank $1$.  Let $\mathscr L_\chi$ be the local system of rank $1$ on $T$ corresponding to $\chi$ as introduced in Remark \ref{cst}. We denote the inclusion $T_{\jmath,\iota} \hookrightarrow T$ by $i$. We apply Lemma \ref{sfc} to  $T_{\jmath,\iota}$ and $i^* \mathscr L_\chi$ to get (4). The functions in (5) are of trace type by our last assumption on $d_T$ (see Definition \ref{dtp}).
\end{proof}

\subsection{Dimension estimation}
\label{dimest}

For $g \in G$, let $\mathcal X_g:=\{cH\in G/H:c^{-1}gc \in H\}$ be the fixed point  subscheme of the spherical variety $\mathcal X=G/H$ under the automorphism $g$. We see $\mathcal X_g \simeq \mathrm T_\circ(g,H)/H$.
 Suppose that $g=su$ is the Jordan decomposition. Then we have $\mathcal X_g=\mathcal X_s \cap \mathcal X_u$ as subschemes of $\mathcal X$. Since $\mathcal X $ is smooth, the subscheme $\mathcal X_s$ fixed by the semisimple element $s$ is smooth (\select{c.f.}, Theorem 13.1 of \cite{Mi}). For $(\jmath,\iota)\in 2^{\mathcal J(T)}\times I(T)$ and $s_1,s_2 \in T_{\jmath,\iota}$, we see that the $\mathrm k$-points of $\mathcal X_{s_1}$ and $\mathcal X_{s_2} $ coincide by Proposition \ref{st}. Hence $\mathcal X_{s_1}=\mathcal X_{s_2}$, and we denote it by $\mathcal X_{\jmath,\iota}$. We define $\mathcal X_{\jmath,\iota,u}:=\mathcal X_{su}$ for some $s \in T_{\jmath,\iota}$ and unipotent $u\in G_\iota$.

 
We write $[\mathcal X \times \mathcal B]_G$ for the subscheme of $G\times \mathcal X \times \mathcal B$ given by $\{(g,xH,rB)\in G\times \mathcal X \times \mathcal B : x^{-1}g x \in H, r^{-1}gr \in B\}$.

Recall that $H$ is a spherical subgroup of $G$ if and only if $\mathcal X$  has only finitely many $B$-orbits (see 2.6 of \cite{Kno}). This is equivalent to $\dim [\mathcal X \times \mathcal B]_G=\dim G$.

We have the following lemma:
\begin{lem} \label{de}
For any $(\jmath,\iota)\in 2^{\mathcal J(T)} \times I(T)$ and unipotent $u \in G_\iota(\mathrm k)$, we have the following dimension estimation
$$\dim  \mathcal X_{\jmath,\iota,u}+\dim  \mathcal B_{\iota,u} +\dim T_{\jmath,\iota} \leq \dim C_{G_\iota}(u)=\dim C_{G}(su),
$$ where $s \in T_\iota(\mathrm k)$.
\end{lem}
\begin{proof}
The last equation is clear.
Define a map $\mathrm c: G \times T_{\jmath,\iota} \times \mathcal X_{\jmath,\iota,u}\times \mathcal B_{\iota,u} \to [\mathcal X \times \mathcal B]_G$ by
$$(g,t,x,r) \mapsto  (g^{-1}tug,g^{-1}x,g^{-1}r).
$$
Set $\mathcal F_{p_0}$ to be the fibre  of $\mathrm c$ over a closed point $p_0=(g_0,x_0,r_0)$.
Let $k(g_0)$ be the (finite) number of elements of the form $tu$ for $t \in T_{\jmath,\iota}(\mathrm k)$ such that $tu$ is conjugate to $g_0$.
Then $\mathcal F_{p_0}$ is isomorphic to the disjoint union of  $k(g_0)$ copies of $C_G(g_0)$. 
If $(g_0,x_0,r_0)$ is in the image of  $\mathrm c$, we see $g_0$ is conjugate to some $su$ for $s\in T_{\jmath,\iota}$.
Hence we have
$$\dim ( G \times T_{\jmath,\iota} \times \mathcal X_{\jmath,\iota,u}\times \mathcal B_{\iota,u} )\leq
\dim [\mathcal X \times \mathcal B]_G +\dim C_{G_\iota}(u)=\dim G +\dim C_{G_\iota}(u).
$$ This completes the proof.
\end{proof}

\begin{defn}\label{def-Ujmathiota}
For $(\jmath,\iota)\in 2^{J(T),F}\times I(T)^F $, We denote the set $$\{[u] \in \mathcal U(G_\iota^F):\dim  \mathcal X_{\jmath,\iota,u}+\dim  \mathcal B_{\iota,u} +\dim T_{\jmath,\iota} = \dim C_{G_\iota}(u)\}$$  by $\mathcal U_{\jmath,\iota}$.
\end{defn}

\subsection{Limits as $\mathcal P_{d_T} \ni \nu  \to \infty$}

The following two lemmas are trivial instances of the Grothendieck trace formula.
\begin{lem}\label{cp}
Let $X_0$ be a separated scheme of finite type over $\mathbb F_q$. Let $r=\dim X_0$. Let $X$ be the pullback of $X_0$ to $\mathrm k$. 
 Let $\mathcal P=\{1+(\nu-1)d \}_{\nu \in \mathbb Z_+}$ be the arithmetic progression starting from $1$ with the gap $d\in \mathbb Z_+$.
Assume that the number $d$ satisfies that $F^d$ acts trivially on the set of irreducible components of $X$. Then we have (see Definition \ref{d-mathcalC} for the definition of $\mathcal C\left(\cdot\right)$.)
$$\lim \limits_{ \mathcal P \ni \nu \to \infty}\frac{|X^{F^\nu}|}{q^{\nu r}}=\mathcal C(X).
$$
\end{lem}

\begin{lem}\label{tslm}
We keep the assumption and the notation as in Remark \ref{cst}.
Fix a character $\eta: R^F \to \Qlb^\times $ and its corresponding sheaf $\mathscr L_\eta$.
 Let $U_0 \stackrel{i}{\hookrightarrow} R_0$ be a $r$-dimensional closed subgroup of $R_0$. Let $U$ be the pullback of $U_0$ to $\mathrm k$. Let $\mathcal P=\{1+(\nu-1)d \}_{\nu \in \mathbb Z_+}$. Assume that the gap $d$ satisfies that $F^d$ acts trivially on the set of components  of $U$, and $\frac{d}{\mathcal C(U)}$ is an integer.
 For $\nu \in \mathcal P$, define 
$$M(\nu)=\sum_{p\in U^{F^\nu}} \mathrm{Tr}(F^\nu,(i^*\mathscr L_\eta)_p).$$
 Then we have
\begin{equation}\nonumber
\lim \limits_{\mathcal P \ni \nu \to \infty} \frac{M(\nu)}{q^{\nu r}}=
\left\{
\begin{aligned}
&0& \text{if the~restricition~of~$\eta$~to ~$U^F$~ is~ nontrivial.}\\
&\mathcal C(U) &\text{otherwise.}
\end{aligned}
\right.
\end{equation}
\end{lem}

\begin{cor}\label{tlm}
 Fix $(\jmath_0,\iota_0) \in 2^{\mathcal J(T),F} \times I(T)^F$ and $\gamma \in N(\iota_0,T)^F$. Let $d_{\jmath_0,\iota_0}=\max \{\dim T_{\jmath_0,\iota_0},0\}$.
 Let $\chi:T^F \to \Qlb^\times$ be a character of $T^F$. 
Define
$M(\nu)=\sum\limits_{s\in T_{\jmath_0,\iota_0}^{F^\nu}}\chi\circ \mathrm N^ \nu(\gamma^{-1}s \gamma)$ for $\nu\in \mathcal P_{d_T}$.
Then the following limit 
$$
\lim _{\mathcal P_{d_T}\ni \nu \to \infty} \frac{M(\nu)}{q^{\nu d_{\jmath_0,\iota_0}}}
 $$
 exists.
\end{cor}

\begin{proof}
The case for empty $T_{\jmath_0,\iota_0}$ is clear. We assume that $T_{\jmath_0,\iota_0}$ is nonempty and has dimension $d_{\jmath_0,\iota_0}$ in the remainder of the proof.
We note that $d_{\jmath_0,\iota_0}=\dim {\dot{T}}_{\jmath_0,\iota_0}$ (see Remark \ref{rmtos} and Definition \ref{tcirc}), since $T_{\jmath_0,\iota_0}$ is an open subscheme of the algebraic group ${\dot{T}}_{\jmath_0,\iota_0}$.

 For positive integer $r$ and any function $f: T^{F^r} \to \mathbb C$, we have
$$
\sum_{x\in T^{F^r}_{\jmath_0,\iota_0} }f(x)=
\sum_{\substack{\jmath \in 2^{\mathcal J(T), {F^r}}\\ \jmath_0 \subset \jmath}}(-1)^{o_{F^r}(\jmath -\jmath_0)}\sum_{x\in ({\dot{T}}_{\jmath,\iota_0})^{F^r}-(\bigcup_{\iota> \iota_0} {\dot{T}}_{\jmath,\iota})^{F^r}} f(x),
$$
where $o_{F^r}(S)$ is the number of $F^r$-orbits of a set $S$, and $\iota >\iota_0$ means that $ G_{\iota_0} \subsetneq G_\iota$. And we have
$$
\sum_{x\in ({\dot{T}}_{\jmath,\iota_0})^{F^r}-(\bigcup_{\iota> \iota_0} {\dot{T}}_{\jmath,\iota})^{F^r}} f(x)= \sum_{x\in (\dot{T}_{\jmath,\iota_0})^{F^r}}f(x) + \sum_{\Theta \subset \Theta_{\iota_0}^{F^r}}(-1)^{|\Theta|}\sum_{x\in ({\dot{T}}_{\jmath,\Theta})^{F^r}} f(x),
$$
where $\Theta_{\iota_0}:=\{\iota \in I(T) : \iota >\iota_0\}$. 
Taking $r$ in $\mathcal P_{d_T}$, we see the sets $2^{\mathcal J(T),F^r}$ and $\Theta_{\iota_0}^{F^r}$ remain unchanged. And it remains to prove that the following limits exist for functions $f_\nu: T^{F^\nu}\to \mathbb C$ defined by $f_\nu(t)=\chi\circ \mathrm{N}^\nu(\gamma^{-1}s \gamma)$
$$\lim_{\mathcal P_{d_T}\ni \nu \to \infty} \frac{ \sum_{x\in ({\dot{T}}_{\jmath,\Theta})^{F^r}} f_\nu(x)}{q^{\nu d_{\jmath_0,\iota_0}}}.
$$
This follows from our assumptions on $d_T$ and  Lemma \ref{tslm} by noting that $\dim {\dot{T}}_{\jmath,\Theta} \leq d_{\jmath_0,\iota_0}$ for $\jmath\subset \jmath_0$ and $\Theta \subset \Theta_{\iota_0}$.
\end{proof}

\begin{defn}\label{dt...}
Let $(\jmath,\iota)\in 2^{\mathcal J(T),F}\times I(T)^F$ and $\bar\gamma  \in \bar N(\iota, T)^F$. Let $\chi:T^F \to \Qlb^\times$ be a character. 
We define $t_{\bar \gamma,\chi,\jmath,\iota}$ to be the limit introduced in Corollary \ref{tlm} for $(\jmath_0,\iota_0)=(\jmath,\iota)$ and a representative $\gamma\in N(\iota,T)^F$   of $\bar \gamma$. (Note that $\bar N(\iota,T)^F=G_\iota^F\backslash N(\iota,T)^F$ by Lang's theorem.)
\end{defn}

\subsection{Periods are constant functions of geometric type for $\mathcal P_{d_T}$}
We use Lemma \ref{const} to compute $\mathrm P(\nu,T,\chi)$.

\begin{thm} \label{main}
Fix a $F$-stable torus $T$ of $G$ and a character $\chi:T^F\to \Qlb^\times$. For $\nu\in \mathcal P_{d_T}$, the function $\mathrm P(\nu, T,\chi)$ is a constant:
$$\mathrm P(\nu,T,\chi)=\sum_{\substack{\jmath \in 2^{\mathcal J(T),F} \\ \iota \in I(T)^{F}}}\sum_{\bar \gamma \in \bar N(\iota,T)^{F}}\sum_{[u]\in \mathcal U_{\jmath,\iota}}
\frac{\mathcal C(\mathrm T(\jmath,\iota,u))\cdot q_{\iota,u,\gamma}\cdot t_{\bar\gamma,\chi,\jmath,\iota}}{|\bar N(\iota,T)^F|\cdot \mathcal C(C_{G_\iota}(u))}.
$$
In particular, we have
$$\langle R_{T,\chi}^G,1_{H^F} \rangle_{H^F}= \mathrm P(1,T,\chi)=
\sum_{\substack{\jmath \in 2^{\mathcal J(T),F} \\ \iota \in I(T)^{F}}}\sum_{\bar \gamma \in \bar N(\iota,T)^{F}}\sum_{[u]\in \mathcal U_{\jmath,\iota}}
\frac{\mathcal C(\mathrm T(\jmath,\iota,u))\cdot q_{\iota,u,\gamma}\cdot t_{\bar\gamma,\chi,\jmath,\iota}}{|\bar N(\iota,T)^F|\cdot \mathcal C(C_{G_\iota}(u))}.
$$
\end{thm}
\begin{proof}
By Lemma \ref{cp}, Corollary \ref{tlm} and Lemma \ref{de}, the following
$$\frac{|\left(\mathrm T(\jmath,\iota)\cap \mathrm T(u,H)\right)^{F^\nu}|\cdot   Q_{\gamma T\gamma^{-1}}^{G_\iota,\nu}(u)}{|\bar N(\iota,T)^{F}|\cdot|C_{G_\iota}(u)^{F^\nu}|\cdot|H^{F^\nu}|}\sum_{s\in T_{\jmath,\iota}^{F^\nu}}\chi\circ \mathrm N^ \nu(\gamma^{-1}s \gamma)$$
as a function of $\nu \in \mathcal P_{d_T}$ has a limit as $\mathcal P_{d_T} \ni \nu \to \infty$, and has a  nonzero limit as $\mathcal P_{d_T}\ni \nu \to \infty$
only if $[u] \in \mathcal U_{\jmath,\iota}$. By Lemma \ref{cp} and Corollary \ref{tlm} and Lemma \ref{meq}, we see
$$\lim _{\mathcal P_{d_T} \ni \nu \to \infty} \mathrm P(\nu,T,\chi) =\sum_{\substack{\jmath \in 2^{\mathcal J(T),F} \\ \iota \in I(T)^{F}}}\sum_{\bar \gamma \in \bar N(\iota,T)^{F}}\sum_{[u]\in \mathcal U_{\jmath,\iota}}
\frac{\mathcal C(\mathrm T(\jmath,\iota)\cap \mathrm T(u,H))\cdot q_{\iota,u,\gamma}\cdot t_{\bar\gamma,\chi,\jmath,\iota}}{|\bar N(\iota,T)^F|\cdot \mathcal C(C_{G_\iota}(u))}.
$$
It remains to use Proposition \ref{jgt} and Lemma \ref{const}.
\end{proof}

\section{Semisimple Elements}\label{s-semisimpleele}

In this section we deduce a formula Theorem \ref{emmain} under some assumption. This assumption is satisfied by a large variety of interesting examples (see the next section). 

\subsection{A lemma  for $[ \mathcal X\times \mathcal B]_G$}

In this subsection we retain the notation in Section \ref{dimest}.
Recall the subscheme $[\mathcal X\times \mathcal B]_G$ of $G\times \mathcal X\times \mathcal B$ given by $\{(g,xH,rB)\in G\times \mathcal X \times \mathcal B : x^{-1}g x \in H, r^{-1}gr \in B\}$ in Section \ref{dimest}. For $w\in B(\mathrm k)\backslash G(\mathrm k)/H(\mathrm k)$, we  set $[\mathcal X\times \mathcal B]_G^w$ to be the locally closed subscheme of $[\mathcal X\times \mathcal B]_G$ given by $\{(g,xH,rB)\in G\times \mathcal X \times \mathcal B : x^{-1}g x \in H, r^{-1}gr \in B,r^{-1}x\in BwH\}$. 

\begin{lem}\label{eqdims}
    Let the notation be as above. For any $w\in B(\mathrm k)\backslash G(\mathrm k)/H(\mathrm k)$, the scheme $[\mathcal X\times \mathcal B]_G^w$ (as a scheme over $\mathrm k$) is equidimensional of dimension $\dim G$. Moreover, the scheme $[\mathcal X\times B]_G^w$ is the disjoint union of its irreducible components.
\end{lem}
\begin{proof}
   Let $\bar w$ be a representative of $w$. Let $O_w$ be the orbit of $(\bar wH, B)\in \mathcal X\times \mathcal B$ under the action of $G$ given by $g\cdot (xH,yB)\mapsto (gxH,gyB)$. Let $R_{\bar w}=B\cap \bar w H\bar w^{-1}$. We have an isomorphism $O_w \simeq G/R_{\bar w}$. The scheme $[\mathcal X\times \mathcal B]_G^w$ fits into the following cartesian square 
   $$
   \xymatrix{
[\mathcal X\times \mathcal B]_G^w \ar[d]\ar[r]& O_w \ar[d]\\
G\times O_w\ar[r]&O_w\times O_w   
   }
   $$
   where the upper horizontal is the projection, the left vertical is the inclusion, the right vertical is the diagonal map, the lower horizontal is given by $(g,x)\mapsto (gx,x)$. Let $R_{\bar w}^{red}$ be the reduced subgroup of $R_{\bar w}$ with the same underlying topological space. Let $O_{w}^r$ be $G/R_{\bar w}^{red}$. We form a similar cartesian square 
   $$
   \xymatrix{
\mathscr A \ar[d]\ar[r]& O_w^r \ar[d]\\
G\times O_w^r\ar[r]&O_w^r\times O_w^r   
   }
   $$ where the right vertical is the diagonal map, the lower horizontal is given by $(g,x)\mapsto (gx,x)$. It suffices to show the desired properties for the scheme $\mathscr A$ since it is visible that the natural map $\mathscr A\to [\mathcal X\times \mathcal B]_G^w$ induces a homeomorphism of the corresponding topological spaces of $\mathrm k$-points. Now we easily see that $O^r_w$ is a smooth scheme and the map $G\times O^r_w \to O^r_w\times O^r_w$ is smooth whose  fibres are isomorphic to $R_{\bar w}^{red}$. Note that a reduced algebraic group scheme over a perfect field is smooth.   
   Consequently the scheme $\mathscr A$ is smooth and equidimensional of dimension $\dim G$, and $\mathscr A$ is the disjoint union of its irreducible components. And we complete the proof.
\end{proof}

\subsection{Around $q_{\iota,u,\gamma}$ and $\mathcal U_{\jmath,\iota}$}
Roughly speaking, the coefficient $q_{\iota,u,\gamma}$ in Theorem \ref{main} is the leading term of the corresponding Green function. It is desired that all $q_{\iota,u,\gamma}$  in the summand of Theorem \ref{main} equal $\pm 1$. In this subsection, we formulate a condition to ensure it is indeed the case.

\begin{defn}\label{desmul}
    Let $K$ be an affine algebraic group over the algebraically closed field $\mathrm k$. We say that $K$ is essentially of multiplicative type if there exists an open dense (in the Zariski topology) subspace of $K(\mathrm k )$ consisting of semisimple elements.
\end{defn}

Some morphism that we encounter in the  proof of Lemma \ref{de} will play an important role in the subsequent subsections. Here we formally define it.
\begin{defn}\label{dc_} Let $T$ be a maximal torus of $G$.
 For $(\jmath,\iota)\in 2^{\mathcal J(T)}\times I(T)$ and a unipotent element $u\in G_\iota(\mathrm k)$, let $\mathrm c_{\jmath,\iota,u}: G \times T_{\jmath,\iota} \times \mathcal X_{\jmath,\iota,u}\times \mathcal B_{\iota,u} \to [\mathcal X \times \mathcal B]_G$ be the map given by
$$(g,t,x,r) \mapsto  (g^{-1}tug,g^{-1}x,g^{-1}r).
$$
\end{defn}

\begin{rmk}\label{rconsco}
    It is well-known that a reductive group  has only finitely many unipotent conjugacy classes. Let $R_{\jmath,\iota,u}$ be the image of $\mathrm c_{\jmath,\iota,u}$. We see that the collection $$\{R_{\jmath,\iota,u}:\text{$(\jmath,\iota)\in 2^{\mathcal J(T)}\times I(T)$, and $u$ is an unipotent element in $G_\iota(\mathrm k)$}\}$$ forms a {\bf finite} partition of $[\mathcal X\times \mathcal B]_G$ into constructible subsets. 
\end{rmk}

\begin{prop}\label{proesmul}
    Let $G$, $H$ and $B$  be as usual. Keep the notation in Theorem \ref{main}.
    Suppose that  for any $g\in G(\mathrm k)$, the algebraic group $B\cap gHg^{-1}$ is essentially of multiplicative type in the sense of Definition \ref{desmul}. Then for any $(\jmath,\iota)\in 2^{\mathcal J(T),F}\times I(T)^F$ and $[u]\in \mathcal U_{\jmath,\iota}$, we have $u=1$.
\end{prop}

\begin{proof} Fix $(\jmath,\iota)\in 2^{\mathcal J(T),F}\times I(T)^F$ and $[u]\in \mathcal U(G_\iota^F)$.
We can see from the proof of Lemma \ref{de} that we have $[u]\in \mathcal U_{\jmath,\iota}$ if and only if the image of the map $\mathrm c_{\jmath,\iota,u}$ is of dimension $\dim G$. Assume that we have $[u]\in \mathcal U_{\jmath,\iota}$, we see that the image of $\mathrm c_{\jmath,\iota,u}$ contains a nonempty open subset $V$ of $\mathrm k$-points of $[\mathcal X\times \mathcal B]_G$. (As we have $\dim [\mathcal X\times \mathcal B]_G=\dim G$.) Let $p_{23}:[\mathcal X\times \mathcal B]_G\to \mathcal X\times \mathcal B$ be the projection. Take $p=(g,x,r)\in V$. Then $p^{-1}(\{(x,r)\})\cap V$ is a nonempty open subset of $p^{-1}(\{(x,r)\})=xBx^{-1}\cap rHr^{-1}$. By the assumption we see that $xBx^{-1}\cap rHr^{-1}$ is essentially of multiplicative type. Hence $p^{-1}(\{(x,r)\})\cap V$ contains a semisimple element $s$. By the definition of $\mathrm c_{\jmath,\iota,u}$ we see that there is $t\in T_{\jmath,\iota}$ so that $tu$ is conjugate to $s$, yielding $u=1$ given the uniqueness of the Jordan decomposition.
\end{proof}

\begin{rmk}
    Informally, the condition
$$\text{for any $g\in G(\mathrm k)$, the algebraic group $B\cap gHg^{-1}$ is essentially of multiplicative type}
$$
    ensures that the the subset $\{(g,xH,yB)\in [\mathcal X\times \mathcal B]_G(\mathrm k): \text{$g$ is not semisimple}\}$ of $[\mathcal X\times \mathcal B]_G(\mathrm k)$ has dimension $\leq \dim G-1$.    \end{rmk}

We have a kind of inverse of Proposition \ref{proesmul}.

\begin{prop}\label{pjesm}
Let $G,H$ and $B$ be as usual. Let $T$ be a maximal torus of $G$.
    Suppose that for any $(\jmath,\iota)\in 2^{\mathcal J(T)}\times I(T)$ and any unipotent $u\in G_\iota(\mathrm k)$ satisfying $\dim T_{\jmath,\iota}+\dim \mathcal B_{\iota,u}+\dim \mathcal X_{\jmath,\iota,u}=\dim C_{G_\iota}(u)$, we have $u=1$. Then for any $g\in G(\mathrm k)$, the algebraic group $B\cap gHg^{-1}$ is essentially of multiplicative type in the sense of Definition \ref{desmul}. 
\end{prop}

\begin{proof}Take $g\in G(\mathrm k)$.
Assume for the contradiction that any dense open subset of $B\cap gHg^{-1}$ contains non-semisimple elements. 

Since $H$ is a spherical subgroup of $G$, we have $\dim [\mathcal X\times \mathcal B]_G=\dim G$. Let $Z$ be the subscheme of $[\mathcal X\times \mathcal B]_G$ defined by $Z:=\{(g,xH,yB)\in [\mathcal X\times \mathcal B]_G:By^{-1}xH=BgH\}$. By Lemma \ref{eqdims}, we have $\dim Z=\dim G$. Let $Y$ be the subscheme of $\mathcal X\times \mathcal B$ defined by $$Y:=\{(xH,yB)\in \mathcal X\times \mathcal B:By^{-1}xH=BgH\}.$$

We have an obvious map $\mathrm p_{23}:Z\to Y$ given by projection to the later two factors. Each fibre of $\mathrm p_{23}$ is isomorphic to $B\cap gHg^{-1}$. In particular, the map $\mathrm p_{23}$ is equidimensional. By Remark \ref{rconsco} and Lemma \ref{eqdims}, we can take a finite collection $\{R_{\jmath_i,\iota_i,u_i}\}_{1\leq i\leq n}$ (where $n$ is an integer) satisfying:
\begin{itemize}
    \item $\dim R_{\jmath_i,\iota_i,u_i}=\dim G$ for $1\leq i\leq n$;
    \item $\bigcup_{1\leq i\leq n} (R_{\jmath_i,\iota_i,u_i}\cap Z)$ is dense in $Z$.
\end{itemize}
We can see from the proof of Lemma \ref{de} that $\dim R_{\jmath_i,\iota_i,u_i}=\dim G$ for $1\leq i\leq n$  only if $$\dim T_{\jmath_i,\iota_i}+\dim \mathcal B_{\iota_i,u_i}+\dim \mathcal X_{\jmath_i,\iota_i,u_i}=\dim C_{G_{\iota_i}}(u_i).$$
By assumption, we have $u_i=1$ for $1\leq i\leq n$. 
Let $R=\bigcup_{1\leq i\leq n} (R_{\jmath_i,\iota_i,u_i}\cap Z)$.
By Lemma \ref{lfden} below, we can take $p=(x_pH,y_pB)\in Y$ so that $\mathrm p_{23}^{-1}(\{p\})\cap R$ is a dense 
(constructible) subset of $x_pHx_p^{-1}\cap y_pBy_p^{-1}$. Since $x_pHx_p^{-1}\cap y_pBy_p^{-1}$ is conjugate to $B\cap gHg^{-1}$ by construction, we see that $\mathrm p_{23}^{-1}(\{p\})\cap R$ contains a non-semisimple element $k\in x_pHx_p^{-1}\cap y_pBy_p^{-1}$. And $k$ is conjugate to an element $t\in T_{\jmath_i,\iota_i}$ for some $1\leq i\leq n$ by construction. This yields a contradiction.
\end{proof}

\begin{lem}\label{lfden}
    Let $f:X\to S$ be an equidimensional surjective map between finite-type separated schemes over a field. Let $V$ be a constructible subset of $X$. Suppose that for any irreducible component $P$ of $X$ with $\dim P=\dim X$, the generic point of $P$ is contained in $V$. Then there exists $s\in S$ so that $f^{-1}(s)\cap V$ is dense in $f^{-1}(s)$.
\end{lem}
\begin{proof}
    By replacing $V$ by a subset we may assume $V$ is open in $X$. Let $Z$ be the complement of $V$ in $X$ endowed with the reduced scheme structure. Assume for the contradiction that for any $s\in S$, the set $f^{-1}(s)\cap V$ is not dense in $f^{-1}(s)$.

    Let $f_0:Z\to S$ be the restriction of $f$. The assumption in the previous paragraph ensures that $f_0$ is likewise surjective.
    Since $f$ is equidimensional, we see that the fibre of $f_0$ has dimension $\dim X- \dim S$. This in turn implies that $\dim Z=(\dim X-\dim S)+\dim S=\dim X$. This is absurd, since any irreducible component $P$ with $\dim P=\dim X$ has its generic point outside $Z$ by construction.
\end{proof}

\begin{rmk}
    The assumption of Proposition \ref{proesmul} is independent of the choice of $B$. Namely, suppose that $B'$ is another Borel subgroup of $G$, and  for  any $g\in G(\mathrm k)$, the algebraic group $B\cap gHg^{-1}$ is essentially of multiplicative type. Then we easily verify that for any any $g\in G(\mathrm k)$, the algebraic group $B'\cap gHg^{-1}$ is likewise essentially of multiplicative type. 
\end{rmk}

\begin{prop}\label{emtg}
    Let $G$ be a connected reductive group over an algebraically closed field and $T\subset B$ be a Borel pair of $G$. Let $U$ be the unipotent radical of $B$.
    Let $\mathrm d:B\to T$ be the map witnessing $T$ as the reductive quotient of $B$ and providing a section for the inclusion $T\hookrightarrow B$.
    Let $B_w$ be a smooth  subgroup of $B$ so that  $B_w$ is 
 essentially of multiplicative (in the sense of Definition \ref{desmul}). Then there exists $b\in B$ satisfying the following:
    \begin{itemize}
        \item[(i)]  $b B_w b^{-1}$ contains $\mathrm d (B_w)$, where we equip $\mathrm d (B_w)$ with the reduced scheme structure;
        \item[(ii)]  the identity component $\mathrm d(B_w)^\circ$ of $\mathrm d(B_w)$ acts (by conjugation) on $b B_w b^{-1}\cap U$ without any fixed points other than the identity $1$.
        \item[(iii)]  $b B_w b^{-1}$ is the semiproduct of $\mathrm d (B_w)$ and the connected unipotent subgroup $b B_w b^{-1}\cap U$ of $U$.
    \end{itemize}
\end{prop}

\begin{proof}
    Let $B_w^\circ$ be the identity component of $B_w$. We may assume that $T_w^\circ:=\mathrm d(B_w^\circ)$ is a subgroup of $B_w^\circ$ by replacing $B_w$ by $b_1B_wb_1^{-1}$ for some $b_1\in B$. Let $\mathrm X(T)$ and $\mathrm X(T^\circ_w)$ be the character lattices of $T$ and $T^\circ_w$ respectively. 
    Let $\Phi_+$ be the set of positive roots corresponding to $B$. Let $\mathrm v_w:\mathrm X(T)\to \mathrm X(T_w^\circ)$ be the map induced by the inclusion $T_w^\circ \hookrightarrow T$.
    Let $\mathrm G_m$ be the split torus of rank $1$ and let $\mathrm X(\mathrm G_m)$ be its character lattice. We fix an identification $\mathrm X(\mathrm G_m)\simeq \mathbb Z$.
    Let $\lambda_0: \mathrm X(T_w^\circ)\to \mathbb Z\simeq \mathrm X(\mathrm G_m)$ be a map of lattices satisfying the following:
    \begin{itemize}
        \item The set $\lambda_0\circ \mathrm v_w(\Phi_+)$ is contained in the monoid of nonnegative integers;
        \item For $\alpha \in \Phi_+$, if $\lambda_0\circ \mathrm v_w(\alpha)=0$, then we have $\mathrm v_w(\alpha)=0$.
    \end{itemize}
We verify that such a morphism $\lambda_0$ exists.    
We abuse the notation by denoting the map of tori induced by the lattice map $\lambda_0$ again by $\lambda_0: \mathrm G_m \to T_w^\circ$.
Let $\lambda : \mathrm G_m \to B_w$ be the composition of $\lambda_0:\mathrm G_m\to T_w^\circ$ and $T_w^\circ \hookrightarrow B_w$.
Then $\lambda$ defines an action of $\mathrm G_m$ on $B_w$ given by
$$
 \mathrm G_m \times B_w \ni(t,b)\mapsto \lambda(t)b\lambda(t)^{-1},
$$
which we denote by $a_\lambda$. This action possesses the following properties (we identify $\mathrm G_m$ with the open subscheme $\mathbb A^1\backslash \{0\}$ of the affine line $\mathbb A^1$):
\begin{itemize}
    \item It extends (in a unique way) to a map $\mathrm l_\lambda :\mathbb A^1\times B_w\to B_w$;
    \item The fixed subscheme of $B_w$ under this action is $C_{B_w}(T_w^\circ)$;
    \item The restriction  of  $\mathrm l_\lambda$ to $\{0\}\times B_w$ gives  a homomorphism of algebraic group $\mathrm r_\lambda: B_w\to C_{B_w}(T_w^\circ)$.
\end{itemize}
By Theorem 13.33 of \cite{Mi}, the homomorphism $r_\lambda$ witnesses $B_w$ as the semiproduct $\mathrm{Ker}(\mathrm r_\lambda)\rtimes C_{B_w}(T_w^\circ)$. The group $\mathrm  {Ker}(\mathrm r_\lambda)$ is connected since for each element $k\in \mathrm {Ker}(\mathrm r_\lambda)$ we have a map of scheme $\mathrm l_{\lambda,k}:\mathbb A^1\to B_w$ given by $x\mapsto \mathrm l_\lambda(x,k)$, linking $k$ with the identity. 
We easily see that $\mathrm{Ker}(\mathrm r_\lambda)$ is contained in the unipotent radical $U$ of $B$ since the map $\mathrm d\circ \mathrm l_\lambda :\mathbb A^1\times B_w \to T$ is visibly equal to  $\mathrm d\circ pr_2:\mathbb A^1\times B_w \to T$.
Since $B_w$ is assumed to be essentially of multiplicative type in the sense of Definition \ref{desmul}, the group $C_{B_w}(T_w^\circ)$ is likewise essentially of multiplicative type by Lemma \ref{lhem}. Let $\wt B_w$ be the 
subgroup of $C_{B_w}(T_w^\circ)$ fitting into the following cartesian diagram:
$$
\xymatrix{
\wt B_w \ar[r]\ar[d]^{\mathrm d_w}& C_{B_w}(T_w^\circ)\ar[d]^{\mathrm d}\\
T_w^\circ \ar[r]& T
}
$$
where the lower horizontal map is the inclusion, and the upper horizontal map  is an open immersion. Clearly, we have $T_w^\circ \subset \wt B_w$, since we see that $T^\circ_w \subset B_w^\circ\cap C_{B_w}(T_w^\circ)$. Moreover, any element $g\in \wt B_w$ can be written as $g=tu$, where $t\in T_w^\circ$ and the unipotent element $u\in U$ commutes with $T_w^\circ$. This implies in particular that the expression $g=su$ is indeed an Jordan  decomposition. Hence the map $\mathrm d_w$ witnesses $\wt B_w$ as the product (of algebraic groups) $(U\cap \wt B_w) \times T_w^\circ$. Given the fact that $\wt B_w$ is essentially of multiplicative type (as it is an open subgroup of $C_{B_w}(T_w^\circ)$), we see that $U\cap \wt B_w$ is the trivial group and $\mathrm d_w$ is indeed an isomorphism. This in turn implies that the map $\mathrm d: C_{B_w}(T_w^\circ)\to T$ is indeed an inclusion of algebraic groups. 

We have seen that $B_w$ is the semiproduct of the connected unipotent group $\mathrm{Ker}(\mathrm r_w)$ and the subgroup $C_{B_w}(T_w^\circ)$, while the latter can be identified with the subgroup $\mathrm d\left(C_{B_w}(T_w^\circ)\right)$ of $T$ (via the map $d$).
And by construction, the group $T_w^\circ$ acts (by conjugation) on the connected unipotent group $\mathrm {Ker}(\mathrm r_w)$ without fixed point other than the identity $1$.
It remains to show that there exists some $b\in B$ so that $b C_{B_w}(T_w^\circ) b^{-1}$ is exactly $\mathrm d\left(C_{B_w}(T_w^\circ)\right)$ as a subgroup of $B$. 
It is routine, as can be seen from the following.
Let $\hat B_w$ be the subgroup  of $B$ fitting into the following cartesian diagram, where the lower horizontal is the inclusion:
$$
\xymatrix{
\hat B_w \ar[r]\ar[d]&B\ar[d]^{\mathrm d}\\
\mathrm d(C_{B_w}(T_w^\circ))\ar[r]&T
}
$$
We have two sections $s_1,s_2$ for the left vertical map given by 
$$t\mapsto t$$
 and 
 $$t\mapsto  s\in C_{B_w}(T_w^\circ)~\text{so that $\mathrm d(s)=t$}.$$
 Then we see from Theorem 16.27 of \cite{Mi} that $s_1$ and $s_2$ differ by conjugation with an element $u\in U$, as desired.
\end{proof}

\begin{lem}\label{lhem}
    Let $h:G_1\to G_2$ be a surjective map of affine smooth algebraic groups. Suppose that $G_1$ is essentially of multiplicative type in the sense of Definition \ref{desmul}. Then $G_2$ is likewise essentially of multiplicative type.
\end{lem}

\begin{proof}
    We  note that the morphism $h$ is an open map  sending an expression $g=su\in G_1$ of Jordan decomposition to the corresponding Jordan decomposition $h(g)=h(s)h(u)\in G_2$ (\select{c.f.} Theorem 9.18 of \cite{Mi}). Since $G_1$ is assumed to be essentially of multiplicative type, we may take an open dense subset $V$ of $G_1(\mathrm k)$ consisting of semisimple elements. 
    We see that $h(V)$ is an open dense subset of $G_2(\mathrm k)$ consisting of semisimple elements, as desired.
\end{proof}

\begin{rmk}\label{resym}
    Suppose that $\theta$ is an involution of $G_0$.
    The pullback of $\theta$ (we abuse the notation by denoting it by $\theta$) to $\mathrm k$ is an involution of $G$.
    We denote the subgroup of $G$ fixed by $\theta$ by $G^\theta$. Let $H$ be the identity component of $G^\theta$. It is well-known that $H$ is a spherical subgroup of $G$. We can show that the pair $(G,H)$ satisfies the assumption of Proposition \ref{proesmul} provided that the characteristic  $p\neq 2$. See Section \ref{e-sym} for details.
\end{rmk}



\subsection{Irreducible components of $[\mathcal X\times \mathcal B]_G$}\label{sicom}
In this subsection, we elaborate the irreducible components of $[\mathcal X\times \mathcal B]_G$ under the assumption of Proposition \ref{proesmul}. 

Let $T$ be an $F$-stable maximal torus of $G$.
Let $B_T$ be a Borel subgroup of $G$ containing $T$. Note that we do not assume that $B_T$ is $F$-stable.

We identify the flag variety $\mathcal B$ with $G/B_T$ in a canonical way, and in particular we may interpret a point $p\in \mathcal B$ as a coset $p=x_pB_T$. For $w\in B_T(\mathrm k)\backslash G(\mathrm k)/H(\mathrm k)$, we  set $[\mathcal X\times \mathcal B]_G^w$ to be the locally closed subscheme of $[\mathcal X\times \mathcal B]_G$ given by $\{(g,xH,rB_T)\in G\times \mathcal X \times \mathcal B : x^{-1}g x \in H, r^{-1}gr \in B_T,r^{-1}x\in B_TwH\}$. 

Let $\mathrm d_T:B_T\to T$ be the map witnessing $T$ as the reductive quotient and providing a section for the inclusion $T\hookrightarrow B_T$.

\begin{defn}\label{dcfra}
For $w\in B_T(\mathrm k)\backslash G(\mathrm k)/H(\mathrm k)$, we denote the set of irreducible components of $\mathrm d_T(B_T\cap wHw^{-1})$ by $\mathfrak C_w$. For $c\in \mathfrak C_w$, we set $[\mathcal X\times \mathcal B]_G^{w,c}$ to be the reduced closed  subscheme of $[\mathcal X\times \mathcal B]_G^w$ consisting of triples $(g,xH,rB_T)$ so that $\mathrm d_T(r^{-1}gr)\in c$. (Here we endow the scheme $\mathrm d_T(B_T\cap wHw^{-1})$ with the reduced scheme structure.)
\end{defn}

It is easy to see that $[\mathcal X\times \mathcal B]_G^{w,c}$ is open in $[\mathcal X\times\mathcal B]_G^w$ as topological spaces.

\begin{lem}\label{irlem}
    Fix $w\in B_T(\mathrm k)\backslash G(\mathrm k)/H(\mathrm k)$.
    Suppose that for arbitrary $g\in G(\mathrm k)$, the algebraic group $B_T\cap gHg^{-1}$ is essentially of multiplicative type in the sense of Definition \ref{desmul}. Then the scheme $[\mathcal X\times \mathcal B]_G^{w,c}$ for $c\in \mathfrak C_w$ is irreducible.
\end{lem}

\begin{proof}
    Let $\mathrm{pr}_{23}:[\mathcal X\times \mathcal B]_G^{w,c}\to \mathcal X\times \mathcal B$ be the projection. The image of $\mathrm{pr}_{23}$ consists of pairs $(xH,rB_T)$ so that $r^{-1} x \in B_TwH$, and  it is isomorphic to the space $G/(B_T\cap wHw^{-1})$. The fibre (endowed with the reduced scheme structure) of $\mathrm{pr}_{23}$ is isomorphic to $c\times (U_T\cap wHw^{-1}) $ as a reduced scheme, where $U_T$ is the unipotent radical of $B_T$. By Proposition \ref{emtg}, the group $U_T \cap wHw^{-1}$ is connected. Consequently, the map $\mathrm{pr}_{23}$ has irreducible image and irreducible fibres of constant dimension, and $[\mathcal X\times \mathcal B]_G^{w,c}$ has a unique irreducible component of maximal dimension. Since $[\mathcal X\times \mathcal B]_G^{w,c}$ is  closed and open in $[\mathcal X\times \mathcal B]_G^w$ as topological spaces, we conclude present lemma by using Lemma \ref{eqdims}.
\end{proof}

\begin{defn}\label{dfgammatbt}
    Suppose that for arbitrary $g\in G(\mathrm k)$, the algebraic group $B\cap gHg^{-1}$ is essentially of multiplicative type in the sense of Definition \ref{desmul}. We define $$\Gamma_{T,B_T}:=\{[\mathcal X \times \mathcal B]_G^{w,c}:\text{ $w\in B_T(\mathrm k)\backslash G(\mathrm k)/H(\mathrm k)$ and $c\in \mathfrak{C}_w$}\}.$$
\end{defn}

\begin{prop}\label{ircp}
Suppose that for arbitrary $g\in G(\mathrm k)$, the algebraic group $B\cap gHg^{-1}$ is essentially of multiplicative type in the sense of Definition \ref{desmul}.
    Then the  set of subschemes $\Gamma_{T,B_T}$ forms a locally closed partion of $[\mathcal X\times \mathcal B]_G$ into irreducible subschemes of dimension $\dim G$. 
\end{prop}

\begin{proof}
    Combining Lemma \ref{irlem} and Lemma \ref{eqdims}.
\end{proof}

\subsection{The incarnation of multi-indices}\label{secincar}
In this subsection, we rephrase Theorem \ref{main} under the assumption of Proposition \ref{proesmul}, changing the complicated multi-indices of Theorem \ref{main} into accessible terms.

We retain the notation in Theorem \ref{main}, and fix a Borel subgroup $B_T$ of $G$ so that $T\subset B_T$. Let $\mathrm d_T:B_T\to T$ and $[\mathcal X\times \mathcal B]_G^{w,c}$  be as 
defined in Section \ref{sicom}.  

The main task of this subsection is to exhibit bijections between $\Phi_T$, $\Gamma_{T,B_T}$ and $\Omega_T$, which is portrayed as the commutative diagram
$$
\xymatrix{
\Phi_T\ar[rr]^{\mathrm V_{T,B_T}}\ar[d]_{\mathrm M_T}&~&\Gamma_{T,B_T}\ar@{.>}[lld]\\
\Omega_T&~&~
}
$$
where $\Gamma_{T,B_T}$ and $\mathrm V_{T,B_T}$ depend on the choice of a Borel subgroup $B_T$ containing $T$, and $\mathrm M_T$ is compatible with the natural actions of the Frobenius endomorphisms. (See Definition \ref{defphi}, \ref{dfgammatbt} and \ref{defomeg}  for the notation.)
With the bijection $\mathrm M_T$ (introduced in Proposition \ref{bjfs}) at hand, we eventually get Theorem \ref{emmain}.

\begin{defn} \label{defphi}
    We define $\Phi_T$
    to be the set consisting of quadruples $(\jmath,\iota,C,X)$ satisfying the following properties:

    \begin{itemize}
        \item $\jmath\in 2^{\mathcal J(T)},~\iota \in I(T)$, and $C$ is an irreducible component of $T_{\jmath,\iota}$;
        \item $X$ is an irreducible  component of $\mathcal X_{\jmath,\iota,1}$ with $\dim X=\dim \mathcal X_{\jmath,\iota,1}$;
        \item $\dim \mathcal X_{\jmath,\iota,1} +\dim \mathcal B_{\iota,1}+\dim T_{\jmath,\iota}=\dim G_\iota$.
    \end{itemize}
  
    Please see Section \ref{dimest} for the notation concerning the dimensional equation. 
\end{defn}

\begin{rmk}\label{rqcom}
    Recall that the scheme $\mathcal B_{\iota,1}$ is by definition the subscheme of $\mathcal B$ fixed by any $s\in T_{\iota}(\mathrm k)$. We may rephrase Proposition 4.4 of \cite{DL} as follows. (See Definition \ref{d-barnugs} and what follows for the definition of $\bar N(\iota,T)$)
    \begin{itemize}
        \item For each  $\bar \gamma \in \bar N(\iota,T)$, there is an irreducible component $\mathcal Q_{\bar \gamma}^\iota$ (defined over $\mathrm k$) of $\mathcal B_{\iota,1}$ consisting of points of the form $gB_T$, so that $g$ is a representative of $\bar \gamma$.
        \item The scheme $\mathcal B_{\iota,1}$ (as a scheme over $\mathrm k$) is the disjoint union of its irreducible components $\mathcal Q_{\bar \gamma}^\iota$ parameterized by $\bar \gamma \in \bar N(\iota,T)$.
    \end{itemize}
\end{rmk}

\begin{defn}
Recall the map $\mathrm c_{\jmath,\iota,u}$ defined in Definition \ref{dc_}.
For an element $\phi=(\jmath,\iota,C,X)\in \Phi_T$, we define $\mathrm c_\phi:G\times C\times X \times \mathcal Q^\iota_1\to [\mathcal X\times \mathcal B]_G$ to be the map given by
$$(g,t,x,r) \mapsto  (g^{-1}tg,g^{-1}x,g^{-1}r).
$$ We note that $\mathrm c_\phi$ is indeed a restriction of $\mathrm c_{\jmath,\iota,1}$ introduced in Definition \ref{dc_}. 
\end{defn}
Recall the set $\Gamma_{T,B_T}$ defined in Definition \ref{dfgammatbt}, which consists of irreducible locally closed subschemes of $[\mathcal X\times \mathcal B]_G$ of dimension $\dim G$.

\begin{lem}\label{delem2}
     Let $G$, $H$ and $B_T$  be as usual. Keep the notation in Theorem \ref{main}.
    Suppose that  for any $g\in G(\mathrm k)$, the algebraic group $B_T\cap gHg^{-1}$ is essentially of multiplicative type in the sense of Definition \ref{desmul}. Take $(\jmath,\iota)\in 2^{\mathcal J(T)}\times I(T)$ and $\bar \gamma \in \bar N(\iota,T)$. Fix an irreducible component $C$ of $T_{\jmath,\iota}$, an irreducible component $X$ of $\mathcal X_{\jmath,\iota,u}$ and an irreducible component $\mathcal Q$ of $\mathcal B_{\iota,u}$. If the image of the restriction of  $\mathrm c_{\jmath,\iota,u}$ to $G\times C\times X\times \mathcal Q$ has dimension $\dim G$, then we have $u=1$ and $(\jmath,\iota,C,X)\in \Phi_T$. 
\end{lem}

\begin{proof}

A similar argument as in Proposition \ref{proesmul} shows that $u=1$ and $\dim \mathcal X_{\jmath,\iota,1} +\dim \mathcal B_{\iota,1}+\dim T_{\jmath,\iota}=\dim G_\iota$. In the remainder of this proof we will show $\dim X=\dim \mathcal X_{\jmath,\iota,1}$.

Since $u=1$, we see that $\mathcal X_{\jmath,\iota,1}$ is smooth by Theorem 13.1 of \cite{Mi}. Consequently, the scheme $X$ is closed and open in $\mathcal X_{\jmath,\iota,1}$.  Also, by Remark \ref{rqcom}, the subscheme $\mathcal Q$ is closed and open in $\mathcal B_{\iota,1}$.
We see from the proof of Lemma \ref{de} that a fibre of the map $\mathrm c_{\jmath,\iota,1}$ is either empty or equidimensional of dimension $\dim G_\iota$.
These, together with the fact that the image of $G\times C\times X\times \mathcal Q$ under $\mathrm c_{\jmath,\iota,1}$ has dimension $\dim G$, imply that $$\dim G+\dim C+\dim X+\dim \mathcal Q=\dim G+\dim G_\iota,$$ yielding $\dim X=\dim \mathcal X_{\jmath,\iota,1}$. (Note that by Lemma \ref{de}, we have $\dim G+\dim C+\dim \mathcal X_{\jmath,\iota,1}+\dim \mathcal Q\leq \dim G+\dim G_\iota$.)
\end{proof}

\begin{defn}\label{defrt}
    The map $\mathrm r_T:[\mathcal X\times \mathcal B]_G\to T$ is defined by sending $(g,xH,yB_T)\in [\mathcal X\times \mathcal B]_G$ to $\mathrm d_T(y^{-1}gy)$.
\end{defn}

\begin{prop}\label{bijin}
Suppose that  for any $g\in G(\mathrm k)$, the algebraic group $B_T\cap gHg^{-1}$ is essentially of multiplicative type in the sense of Definition \ref{desmul}.    We have a bijection $\mathrm V_{T,B_T}: \Phi_T\to \Gamma_{T,B_T}$ characterized by the following two equivalent  conditions:
    \begin{itemize}
        \item For $\phi\in \Phi_T$, the closure (in $[\mathcal X\times \mathcal B]_G$) of the image of $\mathrm c_\phi$ is the closure of $\mathrm V_{T,B_T}(\phi)$;
        \item  For $\phi\in \Phi_T$, let $\eta_{\phi}$ be the generic point of the source of $\mathrm c_\phi$. Then $\mathrm c_\phi$ sends $\eta_\phi$ to the generic point of $\mathrm V_{T,B_T}(\phi)$.
    \end{itemize}
\end{prop}

\begin{proof}
    The two conditions are obviously equivalent. We first show there is a well-defined map $\mathrm V_{T,B_T}$ satisfying the equivalent conditions. 
    Fix $\phi\in \Phi_T$.
    We observe that the source of $c_\phi$ is an irreducible scheme. Given the definition of $\mathrm c_\phi$, we argue as in Lemma \ref{de} to show that the image of $\mathrm c_\phi$ is of dimension $\dim G$. 
    Hence the image of $\mathrm c_\phi$ is an irreducible constructible subset of $[\mathcal X\times \mathcal B]_G$ of dimension $\dim G$.
    By Proposition \ref{ircp}, we see that there is a unique element $[\mathcal X\times \mathcal B]_G^{w,c}\in \Gamma_{T,B_T}$ for some $w\in B_T(\mathrm k)\backslash G(\mathrm k)/H(\mathrm k)$ and $c\in \mathfrak C_w$ so that the generic point of $[\mathcal X\times \mathcal B]_G^{w,c}$ coincides with $\mathrm c_\phi(\eta_\phi)$. Consequently, we have a well-defined map $\mathrm V_{T,B_T}$ satisfying the above listed two equivalent conditions.

    We now show that the map $\mathrm V_{T,B_T}$ is surjective. We have the following:
    \begin{itemize}
        \item[i)] The  images  of maps $\mathrm c_{\jmath,\iota,u}:G\times T_{\jmath,\iota}\times \mathcal X_{\jmath,\iota,u}\times B_{\iota,u}\to [\mathcal X\times \mathcal B]_G$ form a {\bf finite} partition $\{R_{\jmath,\iota,u}\}$ of $[\mathcal X\times \mathcal B]_G$ by Remark \ref{rconsco}. 
        \item[ii)] Fix $(\jmath,\iota)\in 2^{\mathcal J(T)}\times I(T)$ and let $(\jmath,\iota,C,X)\in \Phi_T$. Take a representative $\gamma$ of $\bar \gamma\in\bar N(\iota,T)$. We may assume that $\gamma$ normalizes $T$. There is $(\jmath',\iota') \in 2^{\mathcal J(T)}\times I(T)$ so that $T_{\jmath',\iota'}=\gamma^{-1} T_{\jmath,\iota}\gamma$.
        The image of the restriction of  $\mathrm c_{\jmath,\iota,1}$ to $G\times C\times X \times \mathcal Q^\iota_{\bar \gamma}$ coincides with the image of $\mathrm c_{\phi'}$, where $\phi'=(\jmath',\iota',\gamma^{-1} C \gamma,\gamma^{-1} X)$. 
    \end{itemize}
    The surjectivity of $\mathrm V_{T,B_T}$ follows from Proposition \ref{ircp}, Lemma \ref{delem2} and the above two items.

    We now show that $\mathrm V_{T,B_T}$ is injective. Suppose that we have $\phi_1=(\jmath_1,\iota_1,C_1,X_1),\phi_2=(\jmath_2,\iota_2,C_2,X_2)\in \Phi_T$ so that there are $p_1=(g_1,t_1,x_1,y_1) \in G\times C_1\times X_1 \times \mathcal Q^{\iota_1}_1$ and $p_2=(g_2,t_2,x_2,y_2) \in G\times C_2\times X_2 \times \mathcal Q^{\iota_2}_1$ satisfying $\mathrm c_{\phi_1}(p_1)=\mathrm c_{\phi_2}(p_2)$. We have to show that $\phi_1=\phi_2$. By Lemma \ref{lemrp}, we have $$t_1=\mathrm r_T \circ\mathrm c_{\phi_1}(p_1)=\mathrm r_T \circ \mathrm c_{\phi_1}(p_2)=t_2.$$
    Since $\{T_{\jmath,\iota}\}_{\jmath,\iota}$ forms a partition of $T$, we see that $(\jmath_1,\iota_1)=(\jmath_2,\iota_2)$ and $C_1=C_2$.
    This in turn implies $g_2g_1^{-1}\in C_G(t_1)$, $g_2g_1^{-1}X_1=X_2$ and $g_2g_1^{-1}\mathcal Q_1^{\iota_1}=\mathcal Q_1^{\iota_1}$.
    Hence $g_2g_1^{-1}\in C_G(t_1)\cap G_{\iota_1}B_T$. By Bruhat decomposition, we have $C_G(t_1)\cap G_{\iota_1}B_T=G_{\iota_1}$, \select{c.f.} Proposition 4.4 of \cite{DL}.
    Since $X_1$ and $\mathcal Q_1^{\iota_1}$ are visibly stable under the action of $G_{\iota_1}$, we complete the proof.
\end{proof}

The following lemma is immediate.

\begin{lem}\label{lemrp}
    Let $\phi=(\jmath,\iota,C,X)\in \Phi_T$. Then the composition 
    $$
    G\times C\times X\times \mathcal Q^\iota_1 \stackrel{\mathrm{c}_\phi}{\to} [\mathcal X\times \mathcal B]_G \stackrel{\mathrm r_T}{\to} T 
    $$
  equals the composition of the projection and the obvious inclusion
  $$
  G\times C\times X\times \mathcal Q^\iota_1\to C\hookrightarrow T.
  $$ 
\end{lem}

\begin{rmk}\label{rmopsc} Let $\phi \in \Phi_T$ and $\mathrm V_{T,B_T}(\phi)=[\mathcal X\times \mathcal B]_G^{w,c}$ for some $w\in B_T(\mathrm k)\backslash G(\mathrm k )/H(\mathrm k)$ and $c\in \mathfrak{ C}_w$.
    By Proposition \ref{bijin}, the scheme $\mathrm{c}_\phi^{-1}([\mathcal X\times \mathcal B]_w^{w,c})$ is a locally closed subscheme of $G\times C\times X\times \mathcal Q_1^\iota$ containing the generic point $\eta_\phi$, indicating that $\mathrm{c}_\phi^{-1}([\mathcal X\times \mathcal B]_w^{w,c})$ is open in $G\times C\times X\times \mathcal Q_1^\iota$.
\end{rmk}

\begin{cor}\label{cordc}
Suppose that  for any $g\in G(\mathrm k)$, the algebraic group $B_T\cap gHg^{-1}$ is essentially of multiplicative type in the sense of Definition \ref{desmul}.  
 Let $\phi=(\jmath,\iota,C,X)\in \Phi_T$. Let $\mathrm V_{T,B_T}(\phi)=[\mathcal X\times \mathcal B]_G^{w,c}$ for $w\in B_T(\mathrm k)\backslash G(\mathrm k)/H(\mathrm k )$ and $c\in \mathfrak{C}_w$. Then the closure of $C$ (in both $G$ and $T$) is $c$. 
\end{cor}
\begin{proof}

    By Proposition \ref{bijin}, the closure of the image of $\mathrm c_\phi$ is the closure of $[\mathcal X\times \mathcal B]_G^{w,c}$. Hence $\mathrm c_\phi$ restricts to a dominant map $\mathrm c^{-1}_\phi([\mathcal X\times \mathcal B]_G^{w,c})\to [\mathcal X\times \mathcal B]_G^{w,c}$. Unraveling the definition, we see that $\mathrm r_T$ restricts to a dominant (indeed, surjective) map $[\mathcal X\times \mathcal B]_G^{w,c}\to c.$ Consequently, we have a dominant map $\mathfrak{q}_\phi:\mathrm c^{-1}_\phi([\mathcal X\times \mathcal B]_G^{w,c})\to c$ by composing the above two. 
    The composition (we denote it by $\mathrm p_\phi$) of the inclusion and the projection$$
    \mathrm c^{-1}_\phi([\mathcal X\times \mathcal B]_G^{w,c})
    \hookrightarrow G\times C\times X\times \mathcal Q_1^\iota\to C
    $$ is dominant, since $G\times C\times X\times \mathcal Q_1^\iota$ is irreducible and $\mathrm c^{-1}_\phi([\mathcal X\times \mathcal B]_G^{w,c})$ is a nonempty open subset of $G\times C\times X\times \mathcal Q_1^\iota$ by Remark \ref{rmopsc}. Let $V_C$ be a nonempty open subset of $C$ contained in the image of $\mathrm p_\phi$. 
    Given Lemma \ref{lemrp},     
    We have a series of dominant inclusions $V_C\subset \mathrm {Im}(\mathfrak{q}_\phi) \subset c$, where we denote the image of $\mathfrak{q}_\phi$ by $\mathrm {Im}(\mathfrak{q}_\phi)$. 
    Since $c$ is closed in $G$ and $C$ is irreducible, we have $\bar C= \bar {V}_C=c$, where we denote the closure of $C$ and $V_C$ by $\bar C$ and $\bar {V}_C$ respectively.

\end{proof}

\begin{rmk}\label{rmrbt}
    By Corollary \ref{cordc} and Proposition \ref{bijin}, the sets $$\{\bar C: \text{there is $\phi\in \Phi_T$ whose 3rd factor is $C$}\}$$ and $\bigcup_w \mathfrak{C}_w$ coincide, where $w$ ranges over $B_T(\mathrm k)\backslash G(\mathrm k)/H(\mathrm k)$. In particular, we see $\bigcup_w \mathfrak C_w$ is indeed independent of the choice of the Borel subgroup $B_T$ of $G$ containing $T$. Also, we see that $\bigcup_w \mathfrak C_w$ is a collection of closed subschemes of $T$, which is stable under the action of the Weyl group of $T$ and the action of Frobenius $F$.
\end{rmk}

\begin{defn}\label{defomeg}Suppose that  for any $g\in G(\mathrm k)$, the algebraic group $B_T\cap gHg^{-1}$ is essentially of multiplicative type in the sense of Definition \ref{desmul}.
    We define $\Omega_T$ to be the set consisting of pairs $(c,X)$, satisfying the following:
    \begin{itemize}
        \item $c\in \mathfrak{ C}_w$ for some $w\in B_T(\mathrm k)\backslash G(\mathrm k)/H(\mathrm k)$, see Definition \ref{dcfra} for the definition of $\mathfrak{C}_w$;
        \item Let $\mathcal X^c:=\{xH\in \mathcal X: x^{-1}cx\subset H\}$ be the subscheme of $\mathcal X$ fixed by the scheme $c$. The scheme $X$ is an irreducible component of $\mathcal X^c$ satisfying $\dim X=\dim \mathcal X^c$.
    \end{itemize}
    We remark here that the set $\Omega_T$ is independent of the choice of the Borel subgroup $B_T$ of $G$ containing $T$ by Remark \ref{rmrbt}. 
    The scheme $\mathcal X^c$ is smooth by Theorem 13.1 of \cite{Mi}.   
    And we have a canonical action of the Frobenius operator $F$ on $\Omega_T$ sending $(c,X)\in \Omega_T$ to $(F(c),F(X))$. 
\end{defn}

\begin{rmk} \label{rmdimomega}Keep the assumption as in Proposition \ref{bijin}.
Let $c\in \mathfrak{C}_w$ for some $w\in  B_T(\mathrm k)\backslash G(\mathrm k)/H(\mathrm k)$.
Let $(\jmath,\iota)\in 2^{\mathcal J(T)}\times I(T)$ be such that $T_{\jmath,\iota}\cap c$ is dense in $c$.
    The latter item of Definition \ref{defomeg} amounts to 
    $\dim X=\dim G_\iota-\dim c-\dim \mathcal B_{\iota,1}$ by Section \ref{dimest} and  Proposition \ref{bijin}. Namely, we see that $$\dim X \leq \dim \mathcal X^c=\dim \mathcal X_{\jmath,\iota,1}\leq \dim G_\iota-\dim c-\dim \mathcal B_{\iota,1}$$ by Section \ref{dimest}. And we see from Proposition \ref{bijin} and Corollary \ref{cordc} that $$\dim \mathcal X^c=\dim \mathcal X_{\jmath,\iota,1} = \dim G_{\iota}-\dim c-\dim \mathcal B_{\iota,1}$$
    by taking some element $\phi=(\jmath,\iota,C,X')\in \Phi_T$ with $\mathrm V_{T,B_T}(\phi)=[\mathcal X\times \mathcal B]_G^{w,c}$ and checking the definition of $\Phi_T$.
\end{rmk}

\begin{prop}\label{bjfs}
Suppose that  for any $g\in G(\mathrm k)$, the algebraic group $B_T\cap gHg^{-1}$ is essentially of multiplicative type in the sense of Definition \ref{desmul}.      We have a map $\mathrm M_T: \Phi_T\to \Omega_T$ sending $(\jmath,\iota,C,X)\in \Phi_T$ to $(\bar C,X)\in \Omega_T$, where we denote the closure of $C$ in $G$ by $\bar C$. Then the map $\mathrm M_T$ is a bijection and  compatible with the action of $F$. 
\end{prop}

\begin{proof}
    Decoding the definition, we find that the assignment $\Phi_T\ni(\jmath,\iota,C,X)\mapsto (\bar C,X)\in \Omega_T$ is well-defined by Corollary \ref{cordc} and Remark \ref{rmdimomega}.

    We show that $\mathrm M_T$ is injective.
    Let $(c,X)\in \Omega_T$. Suppose that $(c,X)=(\bar C,X)$ for some $(\jmath,\iota,C,X)\in\Phi_T$.
    We observe that $c$ is irreducible and $\{T_{\jmath,\iota}\}_{\jmath,\iota}$ forms a locally closed partition of $T$. As a result, there is a unique pair $(\jmath_c,\iota_c)\in 2^{\mathcal J(T)}\times I(T)$ so that $T_{\jmath_c,\iota_c}\cap c$ is dense in $c$. Hence, we  have $(\jmath_c,\iota_c)=(\jmath,\iota)$.
    We see that $c\cap T_{\jmath,\iota}$ is an irreducible closed subscheme of $T_{\jmath,\iota}$ (as it is open in $c$, and $c$ is an irreducible closed subscheme of $G$). We have inclusions 
    $$C\subset c\cap T_{\jmath,\iota} \subset T_{\jmath,\iota}
    $$
    while $C$ is an irreducible component of $T_{\jmath,\iota}$.
    And we conclude  $C=T_{\jmath,\iota}\cap c$.

    We show that $\mathrm M_T$ is surjective. Let $(c,X)\in \Omega_T$.
    Suppose that $c\in \mathfrak{C}_w$ for some $w\in B_T(\mathrm k)\backslash G(\mathrm k)/H(\mathrm k)$.
    By Proposition \ref{bijin} and Corollary \ref{cordc}, we have some $\phi_0=(\jmath,\iota,C,X_0)\in \Phi_T$ so that $\mathrm V_{T,B_T}(\phi_0)=[\mathcal X\times \mathcal B]_G^{w,c}$ and $\bar C=c$. In particular $\dim X_0=\dim \mathcal X^c$ by Definition \ref{defphi}. Hence we have $p=(\jmath,\iota,C,X)\in \Phi_T$ with $\mathrm M_T(p)=(c,X)$.

    It is obvious that the map $\mathrm M_T$ is compatible with the action of $F$. And we complete the proof.
\end{proof}

\begin{prop}\label{charcdense}
    Suppose that  for any $g\in G(\mathrm k)$, the algebraic group $B_T\cap gHg^{-1}$ is essentially of multiplicative type in the sense of Definition \ref{desmul}.   The bijection $\mathrm V_{T,B_T}\circ \mathrm M_T^{-1}$ is characterized by the following property:
    \begin{itemize}
        \item Let $\omega=(c,X)\in \Omega_T$. We have $\mathrm V_{T,B_T}\circ \mathrm M_T^{-1}(\omega)=[\mathcal X\times \mathcal B]_G^{w,c}$ for some $w\in B_T(\mathrm k)\backslash G(\mathrm k)/H(\mathrm k)$. Let $O_w=B_T\cdot w$ be the $B_T$-orbit of $w$ in $\mathcal X=G/H$. Then $X\cap O_w$ is an open dense subscheme of $X$.   
    \end{itemize}
\end{prop}
\begin{proof}
By Proposition \ref{bijin} and Proposition \ref{bjfs}, we see that the map $\mathrm V_{T,B_T}\circ \mathrm M_T^{-1}$ is a bijection.

Fix $(c,X)\in \Omega_T$.
    By Proposition \ref{bjfs}, there is a unique $\phi=(\jmath,\iota,C,X)\in \Phi_T$ satisfying $\mathrm M_T(\phi)=\omega$ with $c=\bar C$. 
    Let $[\mathcal X\times \mathcal B]_G^{w,c'}=\mathrm V_{T,B_T}(\phi)$. By Corollary \ref{cordc}, we see that $c'=\bar C=c$. 
    It remains to show there is a nonempty open subset of $X$ contained in $O_w$. 

    Take a point $(g,z,xH,yB_T)\in \mathrm c_\phi^{-1}([\mathcal X\times \mathcal B]_w^{w,c})$. 
    By the definition of $\mathcal Q_1^\iota$, we may take $y\in G_\iota$.
    Let $\wt X$ be the subscheme of $X$ fitting into the following pull-back square
    $$
    \xymatrix{
\wt X \ar[r]\ar[d] & \mathrm{c}_\phi^{-1}([\mathcal X\times \mathcal B]_G^{w,c})\ar[d]\\
    \{(g,z,yB_T)\}\ar[r]& G\times C\times \mathcal Q_1^\iota 
    }
    $$
    where the lower horizontal is the inclusion and the right vertical is the projection. 
   By Remark \ref{rmopsc}, we see that $\wt X$ is a nonempty open subscheme of $X$. Further, we have $y\wt x\in O_w$ for all $\wt x\in \wt X$.  Since $X$ is stable under the action of $G_\iota$, we see that $y \wt X$ is a nonempty open subscheme of $X$ contained in $O_w$, as desired.

\end{proof}

\begin{rmk}\label{r-connectedcomponent}
Take $T,B_T$ and $\mathrm d_T:B_T\to T$ as usual.
    Assume that for each $g\in G(\mathrm k)$, the group $\mathrm d_T(B_T\cap g Hg^{-1})$ is connected. Then $\Gamma_{T,B_T}$ is naturally in bijection with the set $B_T(\mathrm k)\backslash G(\mathrm k)/H(\mathrm k)$ (see Definition \ref{dfgammatbt}).  Under this assumption, the map $\mathrm V_{T,B_T}\circ \mathrm M_T^{-1}$ in Proposition \ref{charcdense} gives a bijection between $\Omega_T$ and the set $B_T(\mathrm k)\backslash G(\mathrm k)/H(\mathrm k)$.
\end{rmk}

\begin{defn}\label{deftphi}
    Recall that the set $\Omega_T$ is equipped with a natural action of the Frobenius operator $F$. (See Definition \ref{defomeg}.) Let $\omega\in (c,X)\in \Omega_T^F$, \select{i.e.}, the scheme $c$ and the scheme $X$ are $F$-stable. Let $\phi=(\jmath,\iota,C,X)\in \Phi_T$ so that $\mathrm M_T(\phi)=\omega$. By Proposition \ref{bjfs}, we have $(\jmath,\iota)\in 2^{\mathcal J(T),F}\times I(T)^F$.
    
    \begin{itemize}
        \item [i)]We define $\sigma_\omega=\sigma_\phi:=\sigma(G_\iota)=\sigma(C_{G}(c)^\circ)$;
        \item [ii)]  Let $\chi:T^F\to \Qlb^\times$ be a character of $T^F$. We define $t_{\omega,\chi}=t_{\phi,\chi}:=\frac{1}{|c^F|}\sum\limits_{t\in c^F}\chi(t)$.
    \end{itemize}
    
\end{defn}

We remark that in Definition \ref{deftphi}, the number $t_{\omega,\chi}$ is either $0$ or a root of unity.
We can now state a refined version of Theorem \ref{main} under the assumption of Proposition \ref{proesmul}.

\begin{thm}\label{emmain}
     Let $G$, $H$ and $B$  be as usual. Keep the notation in Theorem \ref{main}.
    Suppose that  for any $g\in G(\mathrm k)$, the algebraic group $B\cap gHg^{-1}$ is essentially of multiplicative type in the sense of Definition \ref{desmul}. Then for $\nu \in \mathcal P_{d_T}$, the following function is a constant:
    $$
    \mathrm P(\nu,T,\chi)=\sum_{\omega\in \Omega_T^F}(-1)^{\sigma_\omega+\sigma(T)}\cdot t_{\omega,\chi}.
    $$
    In particular, we have $$\langle R_{T,\chi}^G ,1_{H^F}\rangle_{H^F}=\mathrm P(1,T,\chi)=\sum_{\omega\in \Omega_T^F}(-1)^{\sigma_\omega+\sigma(T)}\cdot t_{\omega,\chi}.$$
\end{thm}

\begin{proof}
    Fix $\nu \in \mathcal P_{d_T}$. By Theorem \ref{main}, Proposition \ref{proesmul} and Remark \ref{spc}, we have
    $$\mathrm P(\nu,T,\chi)=\sum_{\substack{\jmath \in 2^{\mathcal J(T),F} \\\iota \in I(T)^{F}\\   [1] \in \mathcal U_{\jmath,\iota} }}
    \sum_{\bar \gamma \in \bar N(\iota,T)^{F}}
\frac{\mathcal C(\mathrm T(\jmath,\iota,1))\cdot (-1)^{\sigma(T)+\sigma(G_\iota)}\cdot t_{\bar\gamma,\chi,\jmath,\iota}}{|\bar N(\iota,T)^F|}.
$$
For a triple $(\jmath,\iota,\bar \gamma)$ in the index set of the above summation, we can take a representative $\gamma_1 \in G(\mathrm k)$ of $\bar \gamma$ so that $\gamma_1^{-1}T\gamma_1=T$. Hence there exists $(\gamma(\jmath),\gamma(\iota))\in 2^{\mathcal J(T)}\times I(T)$ so that $T_{\gamma(\jmath),\gamma(\iota)}=\gamma_1^{-1}T_{\jmath,\iota}\gamma_1$. Given that $\bar \gamma \in \bar N(\iota,T)^F$, we see that $T_{\gamma(\jmath),\gamma(\iota)}$ is $F$-stable, indicating $(\gamma(\jmath),\gamma(\iota))\in 2^{\mathcal J(T),F}\times I(T)^F$.
As we can take a representative $\gamma_2\in G^F$ of $\bar \gamma$, we see that that $\mathcal C (\mathrm T(\jmath,\iota,1))=\mathcal C (\mathrm T(\gamma(\jmath),\gamma(\iota),1))$, $\sigma(G_\iota)=\sigma(G_{\gamma(\iota)})$, $t_{\bar\gamma,\chi,\jmath,\iota}=t_{1,\chi,\gamma(\jmath),\gamma(\iota)}$ and $|\bar N(\iota,T)^F|=|\bar N(\gamma(\iota),T)^F|$.
 Collecting the identical terms, we see
 $$\sum_{\substack{\jmath \in 2^{\mathcal J(T),F} \\\iota \in I(T)^{F}\\   [1] \in \mathcal U_{\jmath,\iota} }}
    \sum_{\bar \gamma \in \bar N(\iota,T)^{F}}
\frac{\mathcal C(\mathrm T(\jmath,\iota,1))\cdot (-1)^{\sigma(T)+\sigma(G_\iota)}\cdot t_{\bar\gamma,\chi,\jmath,\iota}}{|\bar N(\iota,T)^F|}=\sum_{\substack{\jmath \in 2^{\mathcal J(T),F} \\\iota \in I(T)^{F}\\   [1] \in \mathcal U_{\jmath,\iota} }}
\mathcal C(\mathrm T(\jmath,\iota,1))\cdot (-1)^{\sigma(T)+\sigma(G_\iota)}\cdot t_{1,\chi,\jmath,\iota}.
$$
Given Lemma \ref{lctphi} below, we see 
$$\sum_{\substack{\jmath \in 2^{\mathcal J(T),F} \\\iota \in I(T)^{F}\\   [1] \in \mathcal U_{\jmath,\iota} }}
\mathcal C(\mathrm T(\jmath,\iota,1))\cdot (-1)^{\sigma(T)+\sigma(G_\iota)}\cdot t_{1,\chi,\jmath,\iota}=\sum_{\substack{\jmath \in 2^{\mathcal J(T),F} \\\iota \in I(T)^{F}\\   [1] \in \mathcal U_{\jmath,\iota} }}\sum_{j\in J_{\jmath,\iota}}
\mathcal C(\mathrm T(\jmath,\iota,1))\cdot (-1)^{\sigma(T)+\sigma(G_\iota)}\cdot t_j.
$$
(See Lemma \ref{lctphi} for the definition of $J_{\jmath,\iota}$ and $t_j$.)
By Lang's theorem and the definition of $\Phi_T$ (see Definition \ref{defphi}), we have (see Definition \ref{deftphi} for $t_{\phi,\chi}$ and $\sigma_\phi$)
$$
\sum_{\substack{\jmath \in 2^{\mathcal J(T),F} \\\iota \in I(T)^{F}\\   [1] \in \mathcal U_{\jmath,\iota} }}\sum_{j\in J_{\jmath,\iota}}
\mathcal C(\mathrm T(\jmath,\iota,1))\cdot (-1)^{\sigma(T)+\sigma(G_\iota)}\cdot t_j=\sum_{\phi\in \Phi_T^F}(-1)^{\sigma(T)+\sigma_\phi}\cdot t_{\phi,\chi}.
$$
Given the bijection $\mathrm M_T$ introduced in Proposition \ref{bjfs} and Definition \ref{deftphi}, we have
$$\sum_{\phi\in \Phi_T^F}(-1)^{\sigma(T)+\sigma_\phi}\cdot t_{\phi,\chi}=\sum_{\omega\in \Omega_T^F}(-1)^{\sigma_\omega+\sigma(T)}\cdot t_{\omega,\chi},
$$
 as desired.   
\end{proof}

\begin{cor}
     Let $G$, $H$ and $B$  be as usual. Keep the notation in Theorem \ref{main}.
    Suppose that  for any $g\in G(\mathrm k)$, the algebraic group $B\cap gHg^{-1}$ is essentially of multiplicative type in the sense of Definition \ref{desmul}. Then 
     we have
     $$\langle R_{T,1}^G ,1_{H^F}\rangle_{H^F}=\sum_{\omega\in \Omega_T^F}(-1)^{\sigma (T)+\sigma_\omega}.$$
\end{cor}

\begin{lem}\label{lctphi}
    Let $(\jmath,\iota)\in 2^{\mathcal J(T),F}\times I(T)^F$ and $\gamma  \in N(\iota, T)^F$. Let $\chi:T^F \to \Qlb^\times$ be a character. Let $\{C_j\}_{j\in J_{\jmath,\iota}}$ be the set of $F$-stable irreducible components of $T_{\jmath,\iota}$, where $J_{\jmath,\iota}$ is an index set. Let $c_j$ be the closure of $C_j$ in $G$ for $j\in J_{\jmath,\iota}$. Set $$t_j=\frac{1}{|c_j^F|}\sum\limits_{s\in c_j^F}\chi(t)$$ for $j\in J_{\jmath,\iota}$. Then we have $t_{1,\chi,\jmath,\iota}=\sum\limits_{j\in J_{\jmath,\iota}} t_j$. See Definition \ref{dt...} for the definition of $t_{1,\chi,\jmath,\iota}$.
\end{lem}
\begin{proof}
We assume that $T_{\jmath,\iota}$ is a nonempty scheme in the following.
Let $M(\nu)=\sum\limits_{s\in T_{\jmath,\iota}^{F^\nu}}\chi\circ \mathrm N^ \nu(s)$ for $\nu\in \mathcal P_{d_T}$. By Definition \ref{dt...}, we have $$ t_{1,\chi,\jmath,\iota}=
\lim _{\mathcal P_{d_T}\ni \nu \to \infty} \frac{M(\nu)}{q^{\nu d_{\jmath,\iota}}},
 $$
 where $d_{\jmath,\iota}=\dim T_{\jmath,\iota}$. Since $F^{d_T}$ acts trivially on the set of irreducible components of $T_{\jmath,\iota}$, we have $M(\nu)=\sum\limits_{j\in J_{\jmath,\iota}} \sum\limits_{s\in C_j^{F^\nu}}\chi\circ \mathrm N^ \nu(s)$.
 It suffices to show $$\lim_{\mathcal P_{d_T}\ni \nu \to \infty} \frac{1}{q^{\nu d_{\jmath,\iota}}}\sum\limits_{s\in C_j^{F^\nu}}\chi\circ \mathrm N^ \nu(s)=t_j=\frac{1}{|c_j^F|}\sum\limits_{s\in c_j^F}\chi(t)$$ for $j\in J_{\jmath,\iota}.$ In what follows, we fix $j\in J_{\jmath,\iota}$.
 
 Note that $T_{\jmath,\iota}$ is a nonempty (by assumption) open subscheme of an algebraic group of dimension $d_{\jmath,\iota}=\dim T_{\jmath,\iota}$ by Remark \ref{rmtos}. 
Consequently, we have $d_{\jmath,\iota}=\dim c_j=\dim C_j$.
 Since $C_j$ is an open dense subscheme of $c_j$, we see that
  $$\lim_{\mathcal P_{d_T}\ni \nu \to \infty} \frac{1}{q^{\nu d_{\jmath,\iota}}}\sum\limits_{s\in C_j^{F^\nu}}\chi\circ \mathrm N^ \nu(s)=\lim_{\mathcal P_{d_T}\ni \nu \to \infty} \frac{1}{q^{\nu d_{\jmath,\iota}}}\sum\limits_{s\in c_j^{F^\nu}}\chi\circ \mathrm N^ \nu(s)=\lim_{\mathcal P_{d_T}\ni \nu \to \infty} \frac{1}{|c_j^{F^\nu}|}\sum\limits_{s\in c_j^{F^\nu}}\chi\circ \mathrm N^ \nu(s).$$
  By the definition of $d_T$ (Definition \ref{dtp}), we verify that the function $f: \mathcal P_{d_T}\to \Qlb$ given by $$\nu \mapsto \frac{1}{|c_j^{F^\nu}|}\sum\limits_{s\in c_j^{F^\nu}}\chi\circ \mathrm N^ \nu(s) $$ is a constant function taking the value $t_j=\frac{1}{|c_j^F|}\sum\limits_{s\in c_j^F}\chi(t)$. (By Remark \ref{rmtos}, the scheme $c_j$ is an $F$-stable irreducible  component of the algebraic group ${\dot{T}}_{\jmath,\iota}$ introduced in Definition \ref{tcirc}.)
  
\end{proof}

\subsection{Locating $\Omega_T$}\label{seclocateomega}
We retain the notation in the previous subsection. In particular, we fix a (not necessarily $F$-stable) Borel subgroup $B_T$ of $G$ containing 
 the $F$-stable maximal torus $T$.
See Proposition \ref{bijin} and Proposition \ref{bjfs} for the definitions of the bijections $\mathrm V_{T,B_T}$ and $\mathrm M_T$.
Recall that the map $\mathrm d_T: B_T\to T$ provides a section for the inclusion $T\hookrightarrow B$. Let $U_T$ be the unipotent radical of $B_T$. We record the following proposition for future use.
\begin{prop}\label{locateomega}
    Suppose that  for any $g\in G(\mathrm k)$, the algebraic group $B_T\cap gHg^{-1}$ is essentially of multiplicative type in the sense of Definition \ref{desmul}.   Fix $(c,X)\in \Omega_T$ and let $v\in B_T(\mathrm k)\backslash G(\mathrm k)/H(\mathrm k)$ satisfying $[\mathcal X\times \mathcal B]_G^{v,c}=\mathrm V_{T,B_T}\circ \mathrm M_T^{-1}\left((c,X)\right)$. Let $w\in B_T(\mathrm k)vH(\mathrm k)$ so that the algebraic group $B_w:=B_T\cap wHw^{-1}$ contains $T_w:=\mathrm d_T(B_w)$. Then we have the following:
    \begin{itemize}
        \item The point $\bar w=wH\in \mathcal X$ is contained in $X$;
        \item Let $(\jmath,\iota,C,X)=\mathrm M_T^{-1}\left((c,X)\right)\in \Phi_{T}$ and let $B_\iota=B_T\cap G_\iota$. We have a series of dominant inclusions $B_\iota\cdot w H\subset G_\iota \cdot w H \subset X$ of locally closed subschemes of $\mathcal X$. In particular, the group $G_\iota\cap wHw^{-1}$ is a spherical subgroup of $G_\iota$.
    \end{itemize}
\end{prop}

\begin{proof}
    Let $w$ be as in the proposition. All schemes mentioned in this proof are reduced by replacing them with their corresponding reduced closed subschemes with the same underlying spaces.
    The group $B_w$ is assumed to be essentially of multiplicative type. We see  that $U_w:=B_w\cap U_T$ is connected by Proposition \ref{emtg}. And we see from Proposition \ref{emtg} that the identity component $T_w^\circ$ of $T_w$ acts by conjugation on $U_w$ without fixed point other than the identity $1$. 
     And we have $G_\iota=C_G(c)^\circ$ by the definition of $\mathrm M_T$. Note that $c$ is an irreducible component of $T_w$.
    Consequently, we have $G_\iota\cap U_T\cap wHw^{-1}=C_G(c)^\circ \cap U_w \subset C_G(T_w^\circ)^\circ \cap U_w=1$. This in turn implies that  $G_\iota \cap U_w=1$ and $B_\iota\cap wHw^{-1}=G_\iota \cap B_w=T_w$ by Proposition \ref{emtg}. 
    Let $X'$ be the irreducible component of $\mathcal X^c$ containing $\bar w:=wH$ (see Definition \ref{defomeg} for the definition of $\mathcal X^c$, and note that $\mathcal X^c$ is smooth).
    The inclusions $B_\iota\cdot \bar w\subset G_\iota \cdot \bar w\subset X'$ are clear. Note that $T_{\jmath,\iota}$ is equidimensional by Remark \ref{rmtos}, and there is an irreducible component $C$ of $T_{\jmath,\iota}$ satisfying $\bar C=c$ by the definition of $\mathrm M_T$. Hence $\dim T_{\jmath,\iota}=\dim c=\dim T_w$.
    We see that: 
    \begin{itemize}
        \item $\dim B_\iota \cdot \bar w=\dim B_\iota-\dim (B_\iota\cap wHw^{-1})=\dim B_\iota-\dim T_w \leq \dim X'$ by the above argument;
        \item $\dim X'\leq\dim \mathcal X_{\jmath,\iota,1}\leq \dim G_\iota-\dim \mathcal B_{\iota,1}-\dim T_{\jmath,\iota}  =\dim B_\iota -\dim T_w$ by Section \ref{dimest}.
    \end{itemize} 
Consequently, we have $\dim X'=\dim B_\iota-\dim T_{\jmath,\iota}=\dim B_\iota \cdot \bar w$ and $(c,X')\in \Omega_T$ by Remark \ref{rmdimomega}. 
And we deduce that $B_\iota\cdot \bar w$  is an open dense subset of $X'$.
By construction, we have $B_\iota\cdot \bar w\subset B_T \cdot \bar w $. We see from Proposition \ref{charcdense} that $X=X'$. In particular, we have $\bar w\in X$. Moreover, we see that $B_\iota \cdot \bar w$ is dense in the irreducible scheme $X$ and the inclusions $B_\iota \bar w \subset G_\iota \bar w \subset X=X'$ are dominant.
\end{proof}

\section{Examples}\label{s-examp}
In this section we apply Theorem \ref{main} and \ref{emmain} to some examples. We remark here that all of the following examples fit into the picture of Theorem \ref{emmain}.

\subsection{A lemma to compute $t_{\bar\gamma,\chi,\jmath,\iota}$}

The following lemma is a corollary of Lemma \ref{tslm}, which is sometimes useful to calculate $t_{\bar\gamma,\chi,\jmath,\iota}$.

\begin{lem}\label{vtrxlj}
We keep the assumption and symbols as in Remark \ref{cst}. 
Fix a character $\eta: R^F \to \Qlb^\times $ and its corresponding sheaf $\mathscr L_\eta$.
 Let $V_0 \stackrel{j}{\hookrightarrow} R_0$ be a dense open subscheme of $R_0$. Let $V$ be the pullback of $V_0$ to $\mathrm k$.
 For $\nu \in \mathbb Z_+$, define 
$$M(\nu)=\sum_{p\in V^{F^\nu}} \mathrm{Tr}(F^\nu,(i^*\mathscr L_\eta)_p).$$
 Then we have
\begin{equation}\nonumber
\lim \limits_{ \nu \to \infty} \frac{M(\nu)}{q^{\nu \cdot \dim R}}=
\left\{
\begin{aligned}
&0,& \text{if $\eta$  is nontrivial.}\\
&1, &\text{otherwise.}
\end{aligned}
\right.
\end{equation}
\end{lem}
\begin{proof}
Let $Z_0$ be the closed complement of $V_0$ equipped with the reduced scheme structure. Let $Z$ be the pullback of $Z_0$ to $\mathrm k$.
We have
$$M(\nu)=\sum_{p\in V^{F^\nu}} \mathrm{Tr}(F^\nu,(i^*\mathscr L_\eta)_p)=\sum_{p\in R^{F^\nu}} \mathrm{Tr}(F^\nu,(i^*\mathscr L_\eta)_p)-\sum_{p\in Z^{F^\nu}} \mathrm{Tr}(F^\nu,(i^*\mathscr L_\eta)_p).$$
We set $M_0(\nu)=\sum\limits_{p\in Z^{F^\nu}} \mathrm{Tr}(F^\nu,(i^*\mathscr L_\eta)_p)$ for $\nu \in \mathbb Z_+$.
By Lemma \ref{tslm}, it suffices to show 
$$\lim \limits_{ \nu \to \infty} \frac{M_0(\nu)}{q^{\nu \cdot \dim R}}=0.
$$
Since $\dim Z <\dim R$, the above equation is clear.
\end{proof}

\subsection{Parabolic subgroups}\label{e-para}
In this subsection, we set $H$ to be an $F$-stable parabolic subgroup of $G$. Recall that $\mathcal X=G/H$ and $T$ is an $F$-stable maximal torus of $G$. The following proposition should be well-known. We reprove it, using Theorem \ref{emmain}.
\begin{prop}
    Let $G$, $T$, $H$ and $\mathcal X$ be as in the above paragraph. Let $\chi:T^F\to \Qlb^\times$ be a character. Then we have
    $$
    \langle R_{T,\chi}^G,1_{H^F}\rangle_{H^F}=
    \left\{\begin{aligned}
        &0&\text{if $\chi$ is nontrivial,}\\
        &|(\mathcal X^{T})^F|&\text{if $\chi$ is trivial,}
    \end{aligned}\right.
    $$
    where $\mathcal X^T$ is the subscheme of $\mathcal X$ fixed by $T$.
\end{prop}

\begin{proof}
    We can apply Theorem \ref{emmain}, as the assumption there is easy to verify using Bruhat decomposition. Then for each $\omega=(c,X)\in \Omega_T$, we see again by Bruhat decomposition that we have $c=T$. This in turn implies that $X$ is indeed a point of $\mathcal X$.  We see by definition $\sigma_\omega=\sigma(T)$. And the current proposition follows.
\end{proof}

\subsection{The pair $(\mathrm U_4 \times \mathrm U_2,\mathrm U_2\times \mathrm U_2)$} \label{exmain}
Let $\mathrm U_2$ be the unitary group over $\mathbb F_q$ acting on the $2$-dimensional space $V$ over $\mathbb F_{q^2}$ equipped with an Hermitian form $(-,-)_V$. Let $\mathrm U_4$ be the unitary group of $V\oplus V$
 with the Hermitian form $(-,-)_V \oplus (-,-)_V$. 
 Then we have a natural embedding  $i_0:\mathrm U_2 \times \mathrm U_2 \to \mathrm U_4$ of algebraic groups over $\mathbb F_q$. The subgroup  $i:\mathrm U_2 \times \mathrm U_2 \to \mathrm U_4 \times \mathrm U_2$  defined by
$$(g_1,g_2) \mapsto (i_0(g_1,g_2),g_1)
$$
is spherical. As usual, we view the algebraic groups $\mathrm U_4 \times \mathrm U_2$, $\mathrm U_2\times \mathrm U_2$ as algebraic groups over $\mathrm k$ equipped with the corresponding Frobenius endomorphisms.

\begin{rmk}\label{r-swmain}
    Though the pair $(\mathrm U_4 \times \mathrm U_2,\mathrm U_2\times \mathrm U_2)$ fits into the pattern of Theorem \ref{emmain}, we will use Theorem \ref{main} to calculate certain cases of this pair directly.
\end{rmk}

We fix a commutative diagram of schemes over $\mathrm k$ depicted as
$$
\xymatrix{
&\mathrm U_2\times \mathrm U_2 \ar[d]\ar[r]^i &\mathrm U_4 \times \mathrm U_2\ar[d]^{\wt c}\\
&\GL_2\times \GL_2\ar[r]^j &\GL_4\times \GL_2
}
$$
whose vertical arrows are isomorphisms and $j$ is given by $(g_1,g_2) \mapsto (\diag(g_1,g_2),g_1)$. The above diagram does not preserve the Frobenius endomorphism $F$.

\begin{rmk}\label{UUo}
We fix an $F$-stable maximal torus $T=T_{\mathrm U_4}\times T_{\mathrm U_2}$ of $G=\mathrm U_4\times \mathrm U_2$. 

For $(\jmath,\iota) \in 2^{\mathcal J(T),F}\times I(T)^F$ and $[u] \in \mathcal U_{\jmath,\iota}$, we verify  case by case that $u=1$ (details omitted).

For $g\in \GL_4\times \GL_2$, we define the map $c_g:\mathrm U_4\times \mathrm U_2 \to \GL_4\times \GL_2$  by $u \mapsto  g^{-1}\wt c(u)g $.
We may assume that $c_{g_0}(T)$ consists of diagonal matrices for some $g_0 \in\GL_4\times \GL_2$ .
The following are all possible $(\jmath,\iota) \in 2^{\mathcal J(T),F}\times I(T)^F$ such that $\mathcal U _{\jmath,\iota}=\{[1]\}$:
\begin{itemize}
\item[(1)] $c_{g_0}(T_{\jmath,\iota})$ is the $1$-dimensional group whose $\mathrm k$-points are of the form $(\diag(x,x,x,x),\diag(x,x))$;
\item [(2)] $\dim T_{\jmath,\iota}=4$, any element of  $c_{g_0}(T_{\jmath,\iota})$ has $4$ distinct eigenvalues,\\ and  any $s\in c_{g_0}( T_{\jmath,\iota})(\mathrm k)$ is of  the form 
$\left(\diag (x_1,x_2,x_3,x_4),\diag(x_1,x_2)\right)$;
\item [(3)] $T_{\jmath,\iota}$ is $G(\mathrm k)$-conjugate to case (2);
\item[(4)] $\dim T_{\jmath,\iota}=2$, any element of $c_{g_0}(T_{\jmath,\iota})$ has $2$ distinct eigenvalues,\\ and any $s\in c_{g_0}(T_{\jmath,\iota})(\mathrm k)$ is of the form $(\diag(x,y,x,y),\diag(x,y))$;
\item[(5)] $T_{\jmath,\iota}$ is $G(\mathrm k)$-conjugate to case (4);
\item[(6)] $\dim T_{\jmath,\iota}=2$, any element of $c_{g_0}(T_{\jmath,\iota})$  has $2$ distinct eigenvalues,\\ and any $s\in c_{g_0}(T_{\jmath,\iota})(\mathrm k)$ is of the form $(\diag(x,x,x,y),\diag(x,x))$;
\item[(7)] $T_{\jmath,\iota}$ is $G(\mathrm k)$-conjugate to case (6);
\item[(8)] $\dim T_{\jmath,\iota}=3$, any element of $c_{g_0}(T_{\jmath,\iota})$ has $3$ distinct eigenvalues,\\ and any $s\in c_{g_0}(T_{\jmath,\iota})(\mathrm k)$ is of the form $(\diag(x,y,y,z),\diag(x,y))$;
\item[(9)] $T_{\jmath,\iota}$ is $G(\mathrm k)$-conjugate to case (8).
\end{itemize}
\end{rmk}
In the remainder of this subsection, we fix $g_0 \in \GL_4\times \GL_2$ satisfying the assumption in
Remark \ref{UUo}. We will give two corollaries concerning the Deligne-Lusztig characters of $\mathrm U_4\times \mathrm U_2$ induced from some maximal anisotropic torus.

We denote the unitary group of rank $1$ over $\mathbb F_q$ by $\mathrm U_1$. We view $\mathrm U_1$ as an algebraic group over $\mathrm k$ equipped with the Frobenius endomorphism $F: \mathrm U_1 \to \mathrm U_1$ as usual.
We can embed $\prod\limits_{1\leq i \leq 4} \mathrm U_1$ into $\mathrm U_4$ in a natural way, yielding an $F$-stable maximal torus $T_1$ of $\mathrm U_4$. Similarly, we have an $F$-stable  maximal torus $T_2$ of $\mathrm U_2$ that is isomorphic to $\mathrm U_1 \times \mathrm U_1$.
The restriction of $c_{g_0}$ to $T_1\times T_2$ is given by $$\left( a,b   \right) \mapsto (\diag(a_1,a_2,a_3,a_4),\diag(b_1,b_2)),$$ where $a=(a_1,a_2,a_3,a_4)$ along $T_1 \simeq \prod\limits_{1\leq i \leq 4} \mathrm U_1 $ and $b=(b_1,b_2)$ along $T_2 \simeq \mathrm U_1\times \mathrm U_1$.

\begin{cor}
Let $\chi_1,\ldots,\chi_6$ be characters of $\mathrm U_1^F$. Suppose that 
\begin{itemize}
\item[1)]  $\prod \limits_{1 \leq i \leq 6} \chi_1=1$;
\item [2)]$\chi_i \neq 1$ for $1\leq i \leq 4$;
\item [3)]$\chi_i \cdot \chi_j \cdot \chi_k \neq 1$ for $1\leq i<j \leq 4$ and $5\leq k \leq 6$.
\end{itemize}
Let $T=T_1\times T_2$ with $T_1= \prod\limits_{1\leq i\leq 4} \mathrm U_1$ and $T_2 =\mathrm U_1\times \mathrm U_1$, as in the above paragraph.  Let $\eta_1 : T_1^F \to \Qlb^\times$ be the character $\chi_1\boxtimes \chi_2\boxtimes\chi_3 \boxtimes \chi_4$, and $\eta_2 :T_2 \to \Qlb^\times$ be the character $\chi_5\boxtimes \chi_6$. Let $\chi=\eta_1\boxtimes \eta_2$ be the character of $T^F$.
Then we have
$$\langle R_{T,\chi}^{\mathrm U_4\times \mathrm U_2} ,1_{\mathrm U_2^F \times \mathrm U_2^F}\rangle_{\mathrm U_2^F \times \mathrm U_2^F} =-1.
$$
\end{cor}
\begin{proof}
The pair $(\jmath,\iota)  \in 2^{\mathcal J(T),F}\times I(T)^F$ corresponding to (1) in Remark \ref{UUo} contributes to the sum $-1$ in Theorem \ref{main} by Lemma \ref{vtrxlj} according to our assumption 1).
For pairs $(\jmath,\iota)  \in 2^{\mathcal J(T),F}\times I(T)^F$ corresponding to (2), (3), (6), (7), (8), (9) in Remark \ref{UUo}, we see that $t_{\bar\gamma,\chi,\jmath,\iota}=0$ by assumption 2) and Lemma \ref{vtrxlj}. For pairs $(\jmath,\iota)  \in 2^{\mathcal J(T),F}\times I(T)^F$ corresponding to (4) or (5) in Remark \ref{UUo}, we see $t_{\bar\gamma,\chi,\jmath,\iota}=0$ by assumption 3) and Lemma \ref{vtrxlj}. 
\end{proof}

Let $T_3$ be an $F$-stable anisotropic maximal torus of $\mathrm U_4$ which is not isomorphic (as an algebraic group equipped with the Frobenius endomorphism) to $\prod\limits_{1\leq i\leq 4}\mathrm U_1$. Then the torus $T_3$ is isomorphic to $\mathrm U_1 \times R $ for some $R$ being the pullback of an anisotropic torus $R_0$ over $\mathbb F_q$. 
 We have $R^F \simeq \mathrm U_1(\mathbb F_{q^3})$. We have an inclusion $\mathrm e: \mathrm U_1 \to R $ whose restriction to the sets of  $F$-invariant points is the natural inclusion $\mathrm U_1(\mathbb F_q) \hookrightarrow \mathrm U_1(\mathbb F_{q^3})$. The torus $T_3\times T_2$ is an $F$-stable maximal torus of $G=\mathrm U_4\times \mathrm U_2$. The inclusion $\mathrm U_1 \times \mathrm U_1 \simeq Z(\mathrm U_4) \times Z(\mathrm U_2)= Z(G) \hookrightarrow T_3\times T_2$ is given by $(u_1,u_2) \mapsto (a_{u_1},b_{u_2})$, where $a_{u_1}=(u_1,\mathrm e(u_1))$ along $T_3\simeq\mathrm U_1 \times R$ and $b_{u_2}=(u_2,u_2)$ along $T_2 \simeq \mathrm U_1\times \mathrm U_1$.

\begin{cor}\label{exc2}
Keep the notation as introduced in the above paragraph.
Let $\theta_1,\theta_2,\theta_3$ be characters of $\mathrm U_1^F$ and $\theta_4$ be a character of $R^F$.
Let $\eta_4$ be the restriction of the character $\theta_4$ of  $R^F$ to $\mathrm U_1^F$.
Let $\theta_0=\theta_1 \boxtimes \theta_4$ be the character of $T_3^F$. Let $\tau=\theta_2 \boxtimes \theta_3$ be the character of $T_2^F$. Let $S=T_3 \times T_2$ be the $F$-stable torus of $\mathrm U_4\times \mathrm U_2$, as in the above paragraph. Let $\theta =\theta_0\boxtimes \tau$ be the character of $S^F$.  Then we have
$$\langle R_{S,\theta}^{\mathrm U_4\times \mathrm U_2} ,1_{\mathrm U_2^F \times \mathrm U_2^F}\rangle_{\mathrm U_2^F \times \mathrm U_2^F} =
\left\{
\begin{aligned}
&-1,& ~\text{if~$\theta_1$~ is~ nontrivial~and~$\theta_1 \cdot \theta_2 \cdot \theta_3 \cdot \eta_4=1$.}\\
&0, &~\text{otherwise.}
\end{aligned}
\right.
$$
\end{cor}
\begin{proof}
The pair $(\jmath,\iota)  \in 2^{\mathcal J(S),F}\times I(S)^F$ corresponding to (1) in Remark \ref{UUo} contributes to the sum $-1$ (resp. $0$) in Theorem \ref{main}    by Lemma \ref{vtrxlj} if  $\theta_1 \cdot \theta_2 \cdot \theta_3 \cdot \eta_4=1$ (resp. $\theta_1 \cdot \theta_2 \cdot \theta_3 \cdot \eta_4$ is nontrivial). There is no $(\jmath,\iota)  \in 2^{\mathcal J(S),F}\times I(S)^F$ corresponding to (2), (3), (4), (5), (8), (9) in Remark \ref{UUo}. For pairs $(\jmath,\iota)  \in 2^{\mathcal J(T),F}\times I(T)^F$ corresponding to (6) or (7) in Remark \ref{UUo}, the corresponding $t_{\bar\gamma,\theta,\jmath,\iota}=0$ (resp. $t_{\bar\gamma,\theta,\jmath,\iota}=1$) if either $\theta_1$ or $ \theta_2 \cdot \theta_3 \cdot \eta_4$ is nontrivial (resp. both $\theta_1$ and $ \theta_2 \cdot \theta_3 \cdot \eta_4$  are trivial) by Lemma \ref{vtrxlj}. And all terms in Theorem \ref{main} corresponding to (6) and (7) in Remark \ref{UUo} sum up to $0$ (resp. $1$) if either $\theta_1$ or $ \theta_2 \cdot \theta_3 \cdot \eta_4$ is nontrivial (resp. both $\theta_1$ and $ \theta_2 \cdot \theta_3 \cdot \eta_4$  are trivial).
\end{proof}

\subsection{Applications to symmetric spaces}\label{e-sym}
In this subsection, we apply Theorem \ref{emmain} to the connected subgroup of $G$ given by the fixed points of an involution, reproving Theorem 3.3 of \cite{Lusym} under the assumption that the
relevant subgroup is connected.

Let us clarify the notation adopted in this subsection. Let $G$ be as usual. Let $\theta$ be an $F$-stable involution of $G$, \select{i.e.}, $\theta$ commutes with $F$ as an endomorphism of $G$.
Let $K=G^\theta$ be the algebraic group of the fixed points of $\theta$. Let $H=K^\circ$ be the identity component of $K$. 
Recall that $\mathcal X=G/H$ as usual.

It is well-known that $H$ is a spherical subgroup of $G$. See \select{e.g.}, Section 4 of \cite{Spr}. We assume the characteristic $p\neq 2$ in this subsection.

\begin{prop}\label{psymem}
    Let $S\subset B_S$ be a Borel pair of $G$. Suppose that $S$ is stable under the involution $\theta$. Then $B_S\cap H$ is the semiproduct of $S_H:=S\cap H$ and a connected unipotent group $V_S$. Moreover, the identity component $S^\circ_H$ of $S_H$ acts on $V_S$ without any fixed points other than the identity $1$. In particular, the group $B_S\cap H$ is essentially of multiplicative type in the sense of Definition \ref{desmul}.
\end{prop}

\begin{proof}
    These all follow  from Proposition 4.8 of \cite{Spr}, which we now spell out.

    Recall that we have $K=G^\theta$ and $H=K^\circ$. Proposition 4.8 of  \cite{Spr} shows that $K\cap B_S$ is a  semiproduct of $S\cap K$ and a connected unipotent group $V_S$. And our group $H\cap B_S$ is a subgroup of $K\cap B_S$ containing the identity component of $K\cap B_S$. As a result, the group $H\cap B_S$ contains $V_S$, and it is the semiproduct of $S_H=H\cap S$ and $V_S$.

   From Section 4.7 of \cite{Spr} and Lemma \ref{lemsymem} below, we see:
    \begin{itemize}
        \item[] $S_H^\circ$ acts on the Lie algebra of $V_S$ with no nonzero fixed vector.
    \end{itemize}
     We deduce that the subgroup $V_S^{S_H^\circ}$ of $V_S$ fixed by $S_H^\circ$ is finite. Consequently, we see that the morphism $S_H\times  V_S\to B_S\cap H$ given by $(s,v)\mapsto v^{-1}sv$ is dominant with the image consisting of semisimple elements. In particular, the group $B_S\cap H$ is essentially of multiplicative type in the sense of Definition \ref{desmul}. And we see that $V_S^{S_H^\circ}$ is the trivial group by Proposition \ref{emtg}.
\end{proof}

\begin{lem}\label{lemsymem}
    Let $S\subset B_S$ be a Borel pair of $G$. Suppose that $S$ is stable under the involution $\theta$. Let $S_{H}^\circ$ be the identity component of $S_{H}:=S\cap H$. Let $\alpha:S\to \mathrm G_m$ be a root so that the composition $S_{H}^\circ \hookrightarrow S\to \mathrm G_m$ is trivial. Then we have $\theta(\alpha)=-\alpha$ as an element in the character lattice $\mathrm X(S)$ of $S$. 
\end{lem}
\begin{proof}
Let $S^\theta$ be the fixed-point subgroup of $S$.
    We note that $S_H^\circ$ is indeed the identity component of  $S^\theta$. Hence the inclusion $S_H^\circ \hookrightarrow S$ gives a surjection
    $
    \mathrm s:\mathrm X(S)\to \mathrm X(S_H^\circ)
    $
    of the corresponding character lattices. Moreover, the surjection $\mathrm s$ witnesses $\mathrm X(S_H^\circ)$ as the maximal $\theta$-invariant quotient lattice of $\mathrm X(S)$. Tensoring with $\mathbb Q$, we see that the kernel of  $\mathrm s \otimes\mathbb Q$ can be identified with the image $\mathrm {Im}_{\mathbb Q}(\theta-1)$ of the linear map  $\theta-1:\mathrm X(S)\otimes \mathbb Q\to \mathrm X(S)\otimes\mathbb Q$ in a natural way. In particular, we see that $\alpha\in \mathrm {Im}_{\mathbb Q}(\theta-1)$, yielding $\theta(\alpha)=-\alpha$.
\end{proof}

\begin{prop}\label{combtorbit}
    Let $S\subset B_S$ be a Borel pair of $G$. Fix a $B_S$-orbit $O$  in $\mathcal X$. There is  a point $\bar x=xH\in O$ so that $x^{-1}Sx$ is stable under the involution $\theta$. Moreover, the collection of such $\bar x$ forms a single $S$-orbit contained in $O$.
\end{prop}

\begin{proof}
    It is an immediate corollary of Theorem 4.2 of \cite{Spr} (and its proof). See also Proposition 1.3 of \cite{Lusym}.
\end{proof}

\begin{cor}\label{corsymem}
    Let $G$ and $H$ be as introduced at the beginning of this subsection. Then  for any Borel subgroup $B$ of $G$, the intersection $B\cap H$ is essentially of multiplicative type in the sense of Definition \ref{desmul}.
\end{cor}
\begin{proof}
    Combine Proposition \ref{psymem} and Proposition \ref{combtorbit}.
\end{proof}

\begin{rmk}
    By Corollary \ref{corsymem}, the assumption of Proposition \ref{proesmul} is satisfied. And the results in Section \ref{secincar} and Section \ref{seclocateomega} apply. For simplicity, we will use Corollary \ref{corsymem} tacitly in the remainder of this subsection to invoke the results in Section \ref{secincar} and Section \ref{seclocateomega}.
\end{rmk}

Let $T$ be an $F$-stable maximal torus of $G$. Recall the set $\Omega_T$ introduced in Definition \ref{defomeg}. We fix a (not necessarily $F$-stable) Borel subgroup $B_T$ of $G$ containing $T$. Let $\mathrm d_T:B_T\to T$ be the map providing a section for the inclusion $T\hookrightarrow B_T$.

\begin{defn}\label{detow}
    For $w\in B_T(\mathrm k)\backslash G(\mathrm k)/H(\mathrm k)$, we denote the $B_T$-orbit $B_T\cdot w$ in $\mathcal X$ by $O_w$. We set $\mathcal O_w$ to be the $T$-orbit in $O_w$ consisting of points $\bar x =xH\in \mathcal X$ so that $x^{-1}Tx$ is $\theta$-stable.
    See Proposition \ref{combtorbit}.
\end{defn}

\begin{cor}\label{corsym2}
    Let $\omega=(c,X)\in \Omega_T$.  Let $[\mathcal X\times \mathcal B]_G^{w,c}=\mathrm V_{T,B_T}\circ \mathrm M_T^{-1}(\omega)$ for some $w\in B_T(\mathrm k)\backslash G(\mathrm k)/H(\mathrm k)$. 
    Let $\mathrm M_T^{-1}(\omega)=(\jmath,\iota,C,X)\in \Phi_T$.
    Let $\bar x=xH\in \mathcal O_w$.
    Then for any Borel subgroup $B_\iota$ of $G_\iota$ containing $T$, the $B_\iota$-orbit $B_\iota\cdot x$ is an open dense subset of $X$.
\end{cor}
\begin{proof}
Recall that $c$ is an  irreducible component of $\mathrm d_T(B_T\cap xHx^{-1})$. And we have $\mathrm d_T(B_T\cap xHx^{-1})=T\cap xHx^{-1}$ by Proposition \ref{psymem}.

    We first take $B_{\iota,0}=B_T\cap G_\iota$. Then we see from Proposition \ref{locateomega} and Proposition \ref{psymem} (applying to $x^{-1}Tx\subset x^{-1}B_Tx$) that $B_{\iota,0}\cdot x$ is an open dense subset of $X$. And we have $\dim X=\dim B_{\iota,0}-\dim c$.

    Let $B_\iota$ be an arbitrary Borel subgroup of $G_\iota$ containing $T$. We may take a Borel subgroup $B_T'$ of $G$ so that $B_T'\cap G_\iota=B_\iota$. Then Proposition \ref{psymem} applies to the pair $x^{-1}Tx\subset x^{-1}B_T'x$, and we see that $x^{-1}B_\iota x \cap H=x^{-1}Tx\cap H$. In particular, the orbit $B_\iota\cdot x$ contained in $X$ has the dimension $\dim B_\iota-\dim c=\dim B_{\iota,0}-\dim c=\dim X$, yielding the desired conclusion.
\end{proof}

Recall that $\mathfrak{C}_w$ is the collection of components of $\mathrm d_T(B\cap w Hw^{-1})$ introduced in Definition \ref{dcfra}. 

\begin{prop}\label{pinpointsymF}
    The following statements are equivalent for $w\in B_T(\mathrm k)\backslash G(\mathrm k)/H(\mathrm k)$:
    \begin{itemize}
        \item[(i)] $\mathcal O_w$ is $F$-stable;
        \item[(ii)] $\mathcal O_w$ has a point fixed by $F$;
        \item[(iii)]$\mathfrak{C}_w$ is $F$-stable. And for all $F$-stable $c\in \mathfrak{C}_w$, the subscheme $X$ of $\mathcal X$ is $F$-stable, where $(c,X)=\mathrm M_T\circ \mathrm V_{T,B_T}^{-1}([\mathcal X\times \mathcal B]_G^{w,c})$;
        \item[(iv)]There is an $F$-stable $c\in \mathfrak{C}_w$ so that $X$ is $F$-stable, where $(c,X)=\mathrm M_T\circ \mathrm V_{T,B_T}^{-1}([\mathcal X\times \mathcal B]_G^{w,c})$.
    \end{itemize}
\end{prop}
\begin{proof}
    Since $T$ is an $F$-stable connected algebraic group, the equivalence of (i) and (ii) follows from Lang's Theorem.

Assume that we have (ii). Let $\bar y=yH$ be the point mentioned in (ii). We note that $\mathfrak{C}_w$ is the set of components of $\mathrm d_T(B_T\cap xHx^{-1})=T\cap yHy^{-1}$ by Proposition \ref{psymem}, which is $F$-stable.
Let $(c,X)\in \Omega_T$ with $c\in \mathfrak{C}_w$ being $F$-stable. Then we have $(c,F(X))\in \Omega_T$.
Note that $X$ and $F(X)$ are irreducible components of a smooth scheme (which amounts to connected components) by the remark after Definition \ref{defomeg}.
By Corollary \ref{corsym2}, we see that $\bar y\in X\cap F(X)$, indicating $X=F(X)$. And we have shown that (ii) implies (iii).

It is clear that (iii) implies (iv). Note that if $\mathfrak{C}_w$ is $F$-stable, then the identity component of $\mathrm d_T(B_T\cap wHw^{-1})$ is $F$-stable.

    We now show that (iv) implies (ii). Let $\omega=(c,X)=\mathrm M_T\circ \mathrm V_{T,B_T}^{-1}([\mathcal X\times \mathcal B]_G^{w,c})$ with $c$ and $X$ being $F$-stable.
    Let $(\jmath,\iota,C,X)=\mathrm M_T^{-1}(\omega)$.
    Let $\bar x=xH\in \mathcal O_w$.
    Note that $\mathrm M_T$ is a bijection that is compatible with the action of $F$ by Proposition \ref{bjfs}. We see that $(\jmath,\iota)\in 2^{\mathcal J(T),F}\times I(T)^F$.
    Let $B_T'$ be a Borel subgroup of $G$ containing $T$ with $F(B_T')=B_T$.
     We see from Corollary \ref{corsym2} that the $B_T'\cap G_\iota$-orbit generated by $\bar x$ is dense in $X$. Hence the $B_T\cap G_\iota=F(B_T')\cap G_\iota$-orbit generated by $F(\bar x)$ is likewise dense in $X$. Hence $F(x)\in O_w=B_T\cdot w\subset \mathcal X$ by Proposition \ref{charcdense} and the fact $(c,X)=\mathrm M_T\circ \mathrm V_{T,B_T}^{-1}([\mathcal X\times \mathcal B]_G^{w,c})$.  Since $F$ commutes with $\theta$ as an endomorphism of $G$, we see that $F(x)^{-1}T F(x)$ is $\theta$-stable.
      Consequently, we have $F(\bar x)\in \mathcal O_w$ by Proposition \ref{combtorbit} and Definition \ref{detow}.
      Note that $\mathcal O_w$ is a $T$-orbit in $\mathcal X$.
      Let $F(\bar x )=t\bar x$ for $t\in T(\mathrm k)$. By Lang's theorem we have $t=F(s)^{-1} s$ for some $s\in T(\mathrm k)$. We verify that $F(s\bar x)=s\bar x$, exhibiting $s\bar x$ as an $F$-fixed point of $\mathcal O_w$.
\end{proof}

\begin{defn}
    We define a subscheme $\Theta_T$ of $G$ by setting $\Theta_T:=\{g\in G: \theta(g^{-1}Tg)=g^{-1}Tg\}$.
\end{defn}

\begin{rmk}
    We see from Proposition \ref{combtorbit} that $\mathfrak{V}_T:=\Theta_T/H$ is the disjoint union of $\mathcal O_w$ for $w\in B_T(\mathrm k)\backslash G(\mathrm k )/H(\mathrm k)$.
\end{rmk}

\begin{defn}
    Let $\mathcal O$ be a  $T$-orbit in $\mathfrak{V}_T$ with respect to the left action. We define  $\mathfrak{C}_{\mathcal O}$ to be the set of irreducible components of $T\cap xHx^{-1}$ for any $\bar x=xH\in \mathcal O$. 
\end{defn}

We see that any $T$-orbit in $\mathfrak{V}_T$ is indeed of the form $\mathcal O_w$ for some $w\in B_T(\mathrm k)\backslash G(\mathrm k )/H(\mathrm k)$. And we have $\mathfrak{C}_{\mathcal O_w}=\mathfrak{C}_w$ (comparing with Definition \ref{defomeg}).

\begin{defn}
    Let $\mathcal O$ be an $F$-stable $T$-orbit in $\mathfrak{V}_T$. Let $c\in \mathfrak{C}_{\mathcal O}^F$, \select{i.e.}, the irreducible subscheme $c$ of $T$ is $F$-stable. 
    \begin{itemize}
        \item We set $\sigma_{\mathcal O,c}:=\sigma(C_G(c)^\circ)$;
        \item  Let $\chi:T^F\to \Qlb^\times$ be a character. We set $t_{\mathcal O,c,\chi}=\frac{1}{|c^F|}\sum\limits_{t\in c^F}\chi(t)$.
    \end{itemize}
    If we have $\mathcal O=\mathcal O_w$ for some $w\in B_T(\mathrm k)\backslash G(\mathrm k )/H(\mathrm k)$, we also denote $t_{\mathcal O,c,\chi}$ by $t_{w,c,\chi}$ and $\sigma_{\mathcal O,c}$ by $\sigma_{w,c}$.
\end{defn}

\begin{rmk}\label{rmtwc}
  Fix $w\in B_T(\mathrm k)\backslash G(\mathrm k )/H(\mathrm k)$. Suppose that $\mathcal O_w$ is $F$-stable. Let $c\in \mathfrak{C}_{\mathcal O}^F$. Let $\omega=\mathrm M_T\circ \mathrm V_{T,B_T}^{-1}([\mathcal X\times \mathcal B]_G^{w,c})$. Unraveling the definition, we have $\sigma_{w,c}=\sigma_\omega$ and $t_{w,c,\chi}=t_{\omega,\chi}$ (comparing with Definition \ref{deftphi}).
\end{rmk}

\begin{thm}\label{symmain}
    Let $G$ be a reductive group over $\mathbb F_q$. Let $\theta$ be an involution of $G$. Let $H$ be the identity component of $G^\theta$. Let $T$ be an $F$-stable maximal torus of $G$ and $\chi:T^F\to \Qlb^\times$ be a character. We have
$$
\langle R_{T,\chi}^G,1_{H^F}\rangle_{H^F}=
\sum_{\mathcal O} \sum_{c\in \mathfrak{C}_{\mathcal O}^F}(-1)^{\sigma_{\mathcal O,c}+\sigma(T)} t_{\mathcal O,c,\chi},
$$
where the outer sum on the right-hand side is taken over the set of $F$-stable $T$-orbits in $\mathfrak{V}_T$.
\end{thm}

\begin{proof}
    Using Proposition \ref{pinpointsymF}, we pinpoint all $F$-stable elements of $\Omega_T$ in terms of $F$-stable $T$-orbits $\mathcal O$ in $\mathfrak{V}_T$ and $F$-stable elements $c$ in $\mathfrak{C}_{\mathcal O}$. Given Remark \ref{rmtwc}, our theorem is an immediate consequence of Theorem \ref{emmain}.
\end{proof}

\begin{defn}\label{d-chio}
Let $\chi:T^F\to \Qlb^\times$.
     Let $\mathcal O$ be an $F$-stable $T$-orbit in $\mathfrak{V}_T$. Take $xH\in \mathcal O$. We denote the restriction of $\chi$ to $\left(\mathrm d_T(B_T\cap xHx^{-1})\right)^F=(T\cap xHx^{-1})^F$ by $\chi_{\mathcal O}$.
\end{defn}

\begin{rmk}\label{r-epsilon}
Keep the notation in Definition \ref{d-chio}.
    Let $c_0$ be the identity component of $\mathrm d_T(B_T\cap xHx^{-1})$. For $t\in \mathrm d_T(B_T\cap xHx^{-1})$, we denote the connected component of $\mathrm d_T(B_T\cap xHx^{-1})$ containing $t$ by $c_t$.
    Take $xH\in \mathcal O$. By  Proposition 2.3 (c) of \cite{Lusym}, the assignment
    $$
   \left(\mathrm d_T(B_T\cap xHx^{-1})\right)^F \ni t ~\mapsto~\sigma(C_G(c_0)^\circ)\sigma(C_G(c_t)^\circ)
    $$
   determines a character, which we denote by $\epsilon_\mathcal O$. We denote $\sigma_{\mathcal O,c_0}=\sigma(C_G(c_0)^\circ)$ by $\sigma_\mathcal O$ for simplicity.
\end{rmk}

\begin{defn}\label{d-otchi}
    Keep the notation in Remark \ref{r-epsilon}. We define the set $\mathfrak{O}_{T,\chi}$ to be the collection of $F$-stable $T$-orbits $\mathcal O$ in $\mathfrak{V}_T$ satisfying $\chi_\mathcal O=\epsilon_\mathcal O$.
\end{defn}

\begin{rmk}Keep the notation in Remark \ref{r-epsilon}.
  In this remark, we will see that  Theorem \ref{symmain} agrees with Theorem 3.3 of \cite{Lusym}. 
  
  Recall that our subgroup $H$ is always assumed to  be connected.  So we should take the group $K$ in  \cite{Lusym} to be our connected subgroup $H$. 
    We have a natural map $\varphi$ from $T^F\backslash \Theta^F_{T}/H^F$ to the set of $F$-stable $T$-orbits in $\mathfrak{V}_T$.
    Then Theorem 3.3 of \cite{Lusym} claims that $$\langle R_{T,\chi}^G,1_{H^F}\rangle_{H^F}=\sum_{f\in \varphi^{-1}(\mathfrak{O}_{T,\chi})} (-1)^{\sigma_{\varphi(f)}+\sigma(T)}.$$
    We see that Theorem \ref{symmain} and Theorem 3.3 of \cite{Lusym} coincide by noting the following facts.
    \begin{itemize}
    \item[(i)] The set $\varphi^{-1}(\mathfrak{O}_{T,\chi})$ is  in bijection with the set of pairs $(\mathcal O,c)$, where $\mathcal O\in \mathfrak{O}_{T,\chi}$ and $c\in \mathfrak{C}_{\mathcal O}^F$. Moreover, the bijection can be taken to respect the map $\varphi$. (This is an instance of Lang's theorem.)
    \item[(ii)] For an $F$-stable $T$-orbit  $\mathcal O$ in $\mathfrak{V}_T$, the summation $\sum\limits_{c\in \mathfrak{C}_{\mathcal O}^F}(-1)^{\sigma_{\mathcal O,c}+\sigma(T)} t_{\mathcal O,c,\chi}$ is nonzero if and only if $\mathcal O$ is in $\mathfrak{O}_{T,\chi}$
    by Remark \ref{r-epsilon} and Definition \ref{d-otchi}.
    \item[(iii)] Fix an $F$-stable $T$-orbit $\mathcal O$ in $\mathfrak{V}_T$. If $\sum\limits_{c\in \mathfrak{C}_{\mathcal O}^F}(-1)^{\sigma_{\mathcal O,c}+\sigma(T)} t_{\mathcal O,c,\chi}$ is nonzero, we have 
    $$\sum\limits_{c\in \mathfrak{C}_{\mathcal O}^F}(-1)^{\sigma_{\mathcal O,c}+\sigma(T)} t_{\mathcal O,c,\chi}=(-1)^{\sigma_{\mathcal O,c_0}+\sigma(T)}|\mathfrak{C}_{\mathcal O}^F|=(-1)^{\sigma_{\mathcal O}+\sigma(T)}|\mathfrak{C}_{\mathcal O}^F|,$$
    where $c_0$ is the identity component of $T\cap xHx^{-1}$ for any $\bar x=xH\in \mathcal O$. (These follow from Remark \ref{r-epsilon} and Definition \ref{d-otchi}.)
    \end{itemize}
\end{rmk}

\subsection{Reductive subgroups}\label{e-Reeder}
In this subsection, we will see that the complexity $0$ case of \cite{R} fits into the picture of Theorem \ref{emmain}.

Here we clarify the notation. Let $G_1\subset G_2$ be a proper inclusion of connected  reductive groups defined over $\mathbb F_q$. As usual, we view them as $\mathrm k$-schemes endowed with the action of the geometric Frobenius.
The group $G_2$ is assumed to be simple in the sense of Section 6 of \cite{R} (hence it is semisimple with connected Coxeter diagram). Let $G=G_1\times G_2$ and $\wt H=G_1$. Let $\wt H\hookrightarrow G$ be the inclusion given by the diagonal embedding. The characteristic $p$ of our ground field $\mathbb F_q$ is assumed to be good for $G_2$ in this subsection. Let $\mathfrak{g_1}$ ($\mathfrak{g_2}$, resp.) be the Lie algebra of $G_1$
($G_2$, resp.). 
We adopt Assumption 1.3 of \cite{R} in this subsection. (For short, we may adopt a stronger assumption from \select{loc.cit.}: The Killing form of $\mathfrak{g_2}$ is nondegenerate on $\mathfrak{g_1}$. This stronger assumption is easy to be verified in Remark \ref{r-besselpgl} and \ref{r-besselu} when invoking Proposition \ref{p-redesmul}.) 

\begin{defn}\label{d-complexity}
    Let $B_1$ be a Borel subgroup of $G_1$. Let $B_2$ be a Borel subgroup of $G_2$. The complexity $\delta$ is defined to be the minimal codimension of a $B_1$-orbit in $G_2/B_2$ (See Section 6 of \cite{R}). 
\end{defn}

\begin{rmk}\label{rcompl} We keep the notation in Definition \ref{d-complexity}.
    Given Proposition 6.2 of \cite{R}, we see that the complexity $\delta$  equals $\dim G_2-\dim B_2-\dim B_1$. Indeed, in Proposition 6.2 of \cite{R}, we see that $B_1$ has an orbit in $G_2/B_2$ with finite stabilizers. 
\end{rmk}

\begin{rmk}\label{rredsph} 
Keep the notation in Remark \ref{rcompl}.
    Suppose that we have $\delta=0$ in Remark \ref{rcompl}. Then $B_1$ has a dense orbit in $G_2/B_2$. Consequently, we see that $B_1$ is a spherical subgroup of $G_2$ via the composition $B_1 \hookrightarrow G_1\hookrightarrow G_2$. In particular, the double coset $B_2(\mathrm k)\backslash G_2(\mathrm k)/ B_1(\mathrm k)$ is finite. 
    Let $B$ in this remark be the Borel subgroup of $G=G_1\times G_2$ given by $B_1\times B_2$.
    It is obvious that we have a bijection between $B_2(\mathrm k)\backslash G_2(\mathrm k)/B_1(\mathrm k)$ and $B(\mathrm k)\backslash G(\mathrm k)/\wt H(\mathrm k)$.  As a result, the  group $\wt H$ is a spherical subgroup of $G$.
\end{rmk}

In the remainder of this subsection, we assume that $\delta=0$ as in Remark \ref{rredsph}. And we denote $H=\wt H$, which should cause no confusion (since $H$ is a spherical subgroup of $G$ by Remark \ref{rredsph}).

\begin{prop}\label{p-redesmul}
    Let $H$ and $G$ be as introduced in the above paragraph. For any Borel subgroup $B'$ of $G$, the algebraic group $B'\cap H$ is essentially of multiplicative type in the sense of Definition \ref{desmul}.
\end{prop}

\begin{proof}
 We denote the diagonal inclusion by $i:H\hookrightarrow G$ in this proof.
    We will verify the conditions of Proposition \ref{pjesm}. 
    
    Let $T=T_1\times T_2\subset G_1\times G_2=G$ be a maximal torus with $T_1\subset T_2$ (via the inclusion $G_1\hookrightarrow G_2$).
    Let $T_1'$ be the subgroup of $H$ given by $T_1\hookrightarrow G_1 \wt {\to} H$.
    Let $(\jmath,\iota)\in 2^{\mathcal J(T)}\times I(T)$ and unipotent $u\in G_\iota(\mathrm k)$ satisfy $\dim T_{\jmath,\iota}+\dim \mathcal B_{\iota,u}+\dim \mathcal X_{\jmath,\iota,u}=\dim C_{G_\iota}(u)$. We want to show $u=1$. 
   Since $\mathcal X_{\jmath,\iota,u}$ is nonempty, we may assume that $T_{\jmath,\iota}$ is a subscheme of $T_1'$ and $u\in C_{H}(s)^\circ$ for $s\in T_{\jmath,\iota}$,  by replacing $T_{\jmath,\iota}$ and $u$ with some $G(\mathrm k)$-conjugates.
(From here, we fix $s\in T_{\jmath,\iota}$ for simplicity.)
   
    Since $H$ embeds into $G=G_1\times G_2$ via the diagonal, we verify that there is an isomorphism $\mathcal X_{\jmath,\iota,u}\simeq C_{G}(su)/C_{H}(su)$. Then the equation $\dim T_{\jmath,\iota}+\dim \mathcal B_{\iota,u}+\dim \mathcal X_{\jmath,\iota,u}=\dim C_{G_\iota}(u)$ reduces to 
    $$\dim T_{\jmath,\iota}+\dim \mathcal B_{\iota,u}-\dim C_{H}(su)=0.
    $$
    It is clear that $\dim \mathcal B_{\iota,u}=\dim \mathcal B_{G_1}^{su} +\dim \mathcal B_{G_2}^{su}$, where $\mathcal B_{G_1}$ ($\mathcal B_{G_2}$, resp.) is the flag variety of $G_1$ ($G_2$, resp.). (Here we regard $su$ as an element of $G_2$ via the composition $H\wt{\to} G_1\hookrightarrow G_2$, and $\mathcal B_{G_i}^{su}$ is the subscheme of $\mathcal B_{G_i}^{su}$ fixed by $su$ for $i=1,2$.) And the dimensional equation simplifies to 
    $$\dim T_{\jmath,\iota}+\dim \mathcal B_{G_1}^{su} +\dim \mathcal B_{G_2}^{su}-\dim C_{H}(su)=0.
    $$
    
    In the remainder of the proof, we will use extensively the results in \cite{R}. 
    Note that ${T}_{\jmath,\iota}$ is a subscheme of  $Z_\iota$, where $Z_\iota:=C_{G_\iota}(G)\cap T_1'$.
    We have 
    $$0=\dim T_{\jmath,\iota}+\dim \mathcal B_{G_1}^{su} +\dim \mathcal B_{G_2}^{su}-\dim C_{H}(su)\leq \dim Z_{\iota}+\dim \mathcal B_{G_1}^{su} +\dim \mathcal B_{G_2}^{su}-\dim C_{H}(su).$$
    By (6.2), (6.3), (6.4) and Lemma 6.5 of \cite{R}, we have 
    $$
   0\leq \dim Z_{\iota}+\dim \mathcal B_{G_1}^{su} +\dim \mathcal B_{G_2}^{su}-\dim C_{H}(su)\leq \delta=0.$$   
By (6.3), Lemma 6.5 and Lemma 6.1 of \cite{R}, we see that $u=1$ as desired.
\end{proof}

\begin{rmk}\label{r-besselpgl}Let $V$ be an $n+1$-dimensional vector space over $\mathbb F_q$ and let $W$ be an $n$-dimensional subspace of $V$. 
    The diagonal embedding $\mathrm {GL}(W)\hookrightarrow\mathrm {GL}(W)\times \mathrm{PGL}(V)$ satisfies the assumption of Proposition \ref{p-redesmul}, where $\mathrm{GL}(W)\to \mathrm {PGL}(W)$ is the composition of obvious maps $\mathrm {GL}(W)\hookrightarrow \mathrm {GL}(V){\to}\mathrm {PGL}(V).$
    (It is elementary to see that this pair satisfies  Assumption 1.3 of \cite{R}.) In particular, we see that this pair satisfies the assumption of Theorem \ref{emmain}.

   Let $i:\mathrm {GL}(W)\hookrightarrow\mathrm {GL}(V)\times \mathrm {GL}(W)$ be the diagonal embedding. We can show that the pair $(G,H)=(\mathrm {GL}(W)\times \mathrm{GL}(V),\mathrm {GL}(W))$ satisfies the assumption of Theorem \ref{emmain}. To see this, we fix a  Borel subgroup $B$ of $\mathrm {GL}(W)\times \mathrm {GL}(V)$. Let $B'$ be the image of $B$ under the natural projection
 $p:\mathrm {GL}(W)\times \mathrm {GL}(V)\to \mathrm {GL}(W)\times \mathrm {PGL}(V)$. We have a commutative diagram with cartesian squares
 $$
 \xymatrix{
 ? \ar[rr]\ar[d] &~& B \ar[r]\ar[d] &B'\ar[d]\\
 \mathrm {GL}(W)\ar[rr]^i&~&\mathrm {GL}(W)\times \mathrm {GL}(V)\ar[r]^p&\mathrm {GL}(W)\times \mathrm {PGL}(V)
 }
 $$
 where all vertical arrows are inclusions of algebraic groups. We see that $?$ can be simultaneously identified with $B\cap \mathrm {GL}(W)$ in $\mathrm {GL}(W)\times \mathrm {GL}(V)$ and $B'\cap \mathrm {GL}(W)$ in $\mathrm {GL}(W)\times \mathrm {PGL}(V)$. And our assertion follows from the previous paragraph. In particular, we may apply Theorem \ref{emmain} to the pair $(G,H)=(\mathrm {GL}(W)\times \mathrm{GL}(V),\mathrm {GL}(W))$. 
\end{rmk}

\begin{rmk}\label{r-besselu}Fix a positive integer $n$.
    As a variant of Remark \ref{r-besselpgl}, we can also apply Theorem \ref{emmain} to the pair $\mathrm U_n\hookrightarrow \mathrm U_n\times \mathrm U_{n+1}$ given by the diagonal embedding.

    We see from Proposition \ref{p-redesmul} that the diagonal embedding $\mathrm {SO}_{2n} \hookrightarrow \mathrm{SO}_{2n}\times \mathrm {SO}_{2n+1}$ fits into the pattern of Theorem \ref{emmain}. We mention that an elegant formula for this pair (in the regular case in the sense of \cite{R}) is deduced in Section 9 of \cite{R}. (And we see from \select{loc.cit.} that this pair satisfies Assumption 1.3 of \cite{R}.)
    
    Putting these together, we see that the  finite-field analogs of (the basic cases of) Bessel models (in the sense of \cite{GGP}) for Deligne-Lusztig characters fit into the framework of Theorem \ref{emmain}.
\end{rmk}

\appendix
\section{}

In this appendix, we reformulate Deligne-Lusztig characters in terms of character sheaves. 
The main reference for this appendix is \cite{L5}.
We will give an alternative proof of Proposition \ref{jgt}. At the end of this appendix, we will give a proposition concerning certain limits involving Green functions. For a morphism $f:X\to Y$ between schemes, the functors $f_*,f^*,f^!,f_!$ are understood in the derived sense. 
\subsection{Sheaves and Deligne-Lustig characters} \label{sadlc}
Let $G_0$ be a reductive group over $\mathbb F_q$. Let $(T,B)$ be a Borel pair of $G$, so that $T$ is $F$-stable. We fix a map $\mathrm d: B \to T$  witnessing $T$ as the reductive quotient of $B$ and providing a section for the natural inclusion $T \hookrightarrow B$.
 Fix a character $\chi:T^F \to \Qlb^\times$ and its corresponding sheaf $\mathscr L_\chi$ over $T_0$ in the sense of Remark \ref{cst}.

Let $T^{reg}$ be the open subscheme of $T$ consisting of points $t$ satisfying that $C_G(t)=T$. Let $G\times^B B$ be the quotient of $G \times B$ by the action of $B$, where the action is given by $b\cdot(g,b')=(gb^{-1},bb'b^{-1})$. Let $G \times^T T^{reg}$ be the quotient of $G \times T^{reg}$ by the action of $T$, where the action is given by $t\cdot(g,t')=(gt^{-1},tt't^{-1})$. 
Note that $G\times ^T T^{reg}$ is the pullback to $\mathrm k$ of a scheme $G_0 \times^{T_0} T_0^{reg}$ over $\mathbb F_q$, which is defined  in a similar way.
The schemes $G\times ^B B$ and $G\times^T T^{reg}$ are smooth over $\mathrm k$. Let $G^{reg}$ be the union of all conjugations of $T^{reg}$.

 We will abuse the notation by denoting the pullback of $\mathscr L_\chi$ to $T$ again by $\mathscr L_\chi$. 
We fix the notation by the  following diagram with a cartesian square:
\begin{equation}\label{dcs}
\xymatrix{
G\times^{T} T^{reg}\ar[r]^j \ar[d]^{\pi^{reg}}\ar@/^2pc/[rrr]^{\mathrm m} &G \times ^{B} B \ar[d]^{\pi}&G\times B \ar[l]_{p}\ar[r]^{\mathrm d\circ pr_2} &T\\
G^{reg}\ar[r]^o&G&~&~\\
}
\end{equation}
where the map $j$ sends $(g,t)$ to $(g,t)$, the map $\mathrm m$ sends $(g,t)$ to $t$, the map $o$ is the open immersion,  the map $\pi$ sends $(g,b)$ to $g b g^{-1}$, the map $\pi^{reg}$ is the restriction of $\pi$ and  $p$ is the natural quotient map.

There is a sheaf $\mathscr K_{T,B,\chi}$ on $G\times^B B$ satisfying that there is an identification $p^* \mathscr K_{T,B,\chi} \simeq (\mathrm d\circ pr_2)^* \mathscr L_\chi$. The sheaf $j^* \mathscr K_{T,B,\chi}$ is isomorphic to $\mathrm m^* \mathscr L_\chi$, which is indeed a pullback of some sheaf over $G_0 \times^{T_0} T^{reg}_0$. Hence the Frobenius $F$ acts on the local system ${\mathrm m}^* \mathscr L_\chi$ naturally, which gives rise to an action of $F$ on $j^* \mathscr K_{T,B,\chi}$. Let $\mathscr R_{T,B,\chi}=\pi_! \mathscr K_{T,B,\chi}$. The derived object $\mathscr R_{T,B,\chi}[\dim G]$ is the intermediate extension of $\pi^{reg}_!  (j^* \mathscr K_{T,B,\chi})[\dim G]$ to $G$, as the map $\pi$ is small and proper with the top stratum $G^{reg}$. Since we identify the Weil sheaf $j^* \mathscr K_{T,B,\chi}$ with ${\mathrm m}^* \mathscr L_\chi$, we regard $\pi^{reg}_!(j^* \mathscr K_{T,B,\chi})$ as a local system which is independent of $B$. We write $\mathscr R_{T,\chi}$ for $\mathscr R_{T,B,\chi}$ henceforth. 

We note that the maps $\pi^{reg}$ and $o$ are indeed pullbacks of some morphisms defined over $\mathbb F_q$. Hence there is an action of the Frobenius $F$ on $\mathscr R_{T,\chi}$ induced by the canonical action of $F$ on $\pi^{reg}_! ({\mathrm m}^* \mathscr L_\chi)$ by the property of the intermediate extension. The following can be derived from Corollary 2.3.2 of \cite{Lau}. See also Theorem 1.14 (a) of  \cite{L5} for $g\in G^F$ being unipotent, which is adequate for our purpose.

Recall the notation $R_{T,\chi}^{G,\nu}$ defined at the end of Section \ref{dl}.

\begin{thm} \label{dls}
The perverse sheaf $\mathscr R_{T,\chi}[\dim G]$ is the intermediate extension of the perverse sheaf $\pi^{reg}_! ({\mathrm m}^* \mathscr L_\chi)[\dim G]$ to $G$. We equip $\mathscr R_{T,\chi}[\dim G]$ with the action of $F$ induced by the canonical action of $F$ on $\pi^{reg}_! ({\mathrm m}^* \mathscr L_\chi)[\dim G]$. 
Then  $\mathscr R_{T,\chi} [\dim G]$ is pure of weight $\leq \dim G$: that is, for  $\nu \in \mathbb Z_+$, $g\in G^{F^\nu}$ and $i \in \mathbb Z$, the eigenvalues of $F^\nu$ on $\mathscr H^i \left((\mathscr R_{T,\chi}[\dim G])_g\right)$ 
have absolute values $\leq  q^{\nu\cdot(i+\dim G)/{2}}$ for all identifications $\Qlb \simeq \mathbb C$. Moreover,
for  $ \nu \in \mathbb Z_+$ and $g\in G^{F^\nu}$, we have
$$
\mathrm{Tr} \left(F^{\nu},(\mathscr R_{T,\chi})_g\right) =R_{T,\chi}^{G,\nu}(g).
$$
\end{thm}

We give an alternative proof of Proposition \ref{jgt}, using Theorem \ref{dls}.
We denote the subscheme of $G \times \mathcal X$ consisting of pairs $(g,x)$ satisfying that $g$ fixes $x$ by $[\mathcal X]_G$.  Let $\tau: [\mathcal X]_G \to G$ be the projection to the first factor. We see that $\tau$ is indeed a pullback of some morphism $\tau_0$ defined over $\mathbb F_q$, and there is a canonical action of $F$ on $\tau_! \Qlb$ satisfying the following for $\nu\in \mathbb Z_+$ and $g\in G^{F^\nu}$
$$ \mathrm{Tr}\left(  F^\nu,(\tau_! \Qlb)_g \right)= \mathrm{Ind}_{H^{F^\nu}}^{G^{F^\nu}}(1_{H^{F^\nu}})(g).
$$
Consequently, we have
$$\mathrm P(\nu,T,\chi) \cdot |G^{F^\nu}|= \mathrm{Tr} (F^\nu,H^\bullet_c(\tau_! \Qlb \otimes \mathscr R_{T,\chi}))
$$
by Theorem \ref{dls} and the Grothendieck trace formula, where $H^\bullet_c(\tau_! \Qlb \otimes \mathscr R_{T,\chi})$ is the virtual module $$\sum\limits_{i \in \mathbb Z}(-1)^i H^i_c (\tau_! \Qlb \otimes \mathscr R_{T,\chi}).$$
Applying Lemma \ref{sfc}, we see that the function $\nu \mapsto |G^{F^\nu}|$ is of trace type in the sense of Definition \ref{gt}. Obviously we have $|G^{F^\nu}|\neq 0$ for $\nu \in \mathbb Z_+$.
This completes the proof of Proposition \ref{jgt} (indeed we prove Proposition \ref{jgt} for the arithmetic progression $\mathbb Z_+$ here).

\begin{rmk}
    A similar argument is adopted in \cite{S}.
\end{rmk}

\begin{rmk}
     The assumption on the size of $q$ (see Section \ref{conv}) aims to guarantee that Theorem \ref{dls} holds.
\end{rmk}

\subsection{Transport the structure}\label{tts.}
Throughout this subsection, we fix a  Borel pair $(\mathrm T_0,\mathrm B_0)$ of $G_0$.
The pullback of $\mathrm T_0$ (resp. $\mathrm B_0$) to $\mathrm k$ is denoted by $\mathrm T$ (resp. $\mathrm B$).
 We denote the Weyl group of $G$ by $\mathrm W_G$. Then
 the Frobenius $F$ acts naturally on $\mathrm W_G$.

 We run the process of the previous subsection for trivial $\chi: \mathrm T^F \to \Qlb^\times$. Note that the Diagram (\ref{dcs}) for $(\mathrm T,\mathrm B )$ is indeed a pullback of a corresponding diagram defined over $\mathbb F_q$. Hence there is a canonical action of $F$ on $\mathscr R_{\mathrm T,1}$, which agrees with the action defined in Section \ref{sadlc} over $G^{reg}$. By the property of intermediate extension, we see that the canonical action of $F$ on $\mathscr R_{\mathrm T,1}$ and the action defined in Section \ref{sadlc} on $\mathscr R_{\mathrm T,1}$ coincide.

 Fix an $F$-stable maximal torus $T$ of $G$ and a (not necessarily $F$-stable) Borel subgroup $B$ of $G$ containing $T$. Suppose that $T=g_{T} \mathrm T g_{T}^{-1}$ for some $g_T \in G(\mathrm k)$. Then $\mathrm l(g_T):=g^{-1}_TF(g_T)$ lies in $N_G(T,T):=\{g\in G:gTg^{-1}=T\}$. And the $F$-conjugacy   class  of the reduction $w(g_T)$ of $g^{-1}_T F(g_T)$ in $\mathrm W_G =\bar N_G(\mathrm T,\mathrm T):=N_G(\mathrm T,\mathrm T)/\mathrm T$ is independent of the choice of $g_T$, which we denote by $[T]=[w(g_T)]$.  Here we adopt the following definition.
 \begin{defn}
     We say that two elements $w_1,w_2\in \mathrm {W}_G$ are in the same $F$-conjugacy class if  there exists $w\in\mathrm W_G$ satisfying $w_1=w^{-1}w_2F(w)$.
 \end{defn}

We have a commutative diagram
$$
\xymatrix
{
G\times^T  T^{reg} \ar[r]^{\mathrm t}\ar[d]^{\pi^{reg}_T}&G\times^{\mathrm T} \mathrm T^{reg}\ar[ld]^{\pi^{reg}_{\mathrm T}}\\
G^{reg}&~
}
$$ where 
$\mathrm t$ is given by sending $(g,t)$ to $(gg_T,g_T^{-1}tg_T )$, and $\pi^{reg}_T$ (resp. $\pi^{reg}_\mathrm T$) is the map $\pi^{reg}$ in the previous subsection corresponding to $T$ (resp. $\mathrm T$). 

The map $\pi^{reg}_{\mathrm T}$ witnesses $G\times^\mathrm T \mathrm T^{reg}$ as a $\mathrm W_G$-torsor over $G^{reg}$, where the $\mathrm W_{G}$-action is given by $w\cdot(g,t)=(g\dot w^{-1},\dot w t\dot w^{-1})$ for any lifting $\dot w \in N_G(\mathrm T,\mathrm T)$ of $w \in \mathrm W_G$. 
 We see that the map $\mathrm t \circ F \circ \mathrm t^{-1}$ sends $(g,t)$ to $(F(g)(\mathrm l(g_T))^{-1},\mathrm l(g_T)F(t) (\mathrm l(g_T))^{-1})$, which is exactly the action of $w(g_T)\circ F$ on the $\mathrm W_G$-torsor $G \times^{\mathrm T} \mathrm T^{reg}$ over $G^{reg}$. 

For $w\in \mathrm W_G$, let $\wt w:\pi_{\mathrm T,!}^{reg}\Qlb\to\pi_{\mathrm T,!}^{reg}\Qlb$ be the map given by the composition 
 $$\pi_{\mathrm T,!}^{reg}\Qlb\stackrel{adj}{ \to}\pi_{\mathrm T,!}^{reg}w_*w^*\Qlb\to\pi_{\mathrm T,!}^{reg}w_*\Qlb\to \pi_{\mathrm T,!}^{reg}\Qlb,$$
 where the first is given by adjunction, the second is given by the natural $w^*\Qlb\simeq \Qlb$, the third is given by  $w_*\simeq w_!$ and $\pi_{\mathrm T}^{reg}\circ w=\pi_{\mathrm T}^{reg}$. We denote the corresponding intermediate extension to $G$ again by $\wt w: \mathscr R_{\mathrm T,1}\to \mathscr R_{\mathrm T,1}$ and likewise for its stalks.

 Combining Theorem \ref{dls}, we have the following.
\begin{prop}\label{ts} Keep the notation as in the previous paragraphs.
Fix unipotent $u \in G^F$.   For any $F$-stable maximal torus $T$ with $[T]=[w]$,  
 we have
$$\mathrm{Tr}\left(F, (\mathscr R_{T,1})_u\right)= \mathrm{Tr}\left( F\circ \wt w,(\mathscr R_{\mathrm T,1})_u   \right),
$$
where the action of $F$ is  defined in Section \ref{sadlc} (hence on the right hand side, the action of $F$ is canonical by the first paragraph of Section \ref{tts.}).
\end{prop}


\subsection{Sheaves and Green functions}
 We keep the notation as in the previous subsection.  Fix a unipotent element $u\in G^F$. Let $\mathcal B^u:=\{x\mathrm B \in G /\mathrm B: x^{-1}ux \in \mathrm B\}$ be the subscheme of the flag variety fixed by $u$. We denote the dimension of $\mathcal B^u$ by $d_u$. Let $n_u^\circ \in \mathbb Z^+$ be a positive integer satisfying that $F^{n_u^\circ}$ acts trivially on the set of irreducible components of $\mathcal B^u$ and that $F^{n_u^\circ}$ acts trivially on $\mathrm W_G$. We set $n_u=n_u^\circ \cdot |\mathrm W_G|$.
Let $\mathcal P_{n_u}:=\{1+(\nu-1)n_u\}_{\nu \in \mathbb Z^+}$ be the arithmetic progression starting from $1$ with the gap $n_u$.
Then we have the following: 
\begin{prop}\label{lmq}
We keep the notation as in the above paragraph.
Fix a $F$-stable maximal  torus $T$ of $G$.  Suppose that $[w]=[T]$ for $w\in \mathrm W_G$.
Then we have 
$$
\lim_{\mathcal P_{n_u} \ni \nu \to \infty} \frac{Q^{G,\nu}_{T}(u)}{q^{\nu d_u}}=\frac{1}{q^{d_u}}\cdot
\mathrm{Tr}\left(F\circ \wt w, \mathscr H^{2d_u}(\mathscr R_{\mathrm T,1})_u \right).
$$
\end{prop}
\begin{proof}
Recall that the action of $F$ on $\mathscr R_{\mathrm T,1}$ is the canonical one (see the first paragraph of Section \ref{tts.}). We may assume $w=w(g_T)$ as in the second paragraph of Section \ref{tts.}.  Since $F^{n_u}$ acts trivially on the set of irreducible components $\mathcal B^u$, we see that $F^{n_u}$ acts as the multiplication by  $q^{n_u d_u}$ on $\mathscr H^{2 d_u}(\mathscr R_{\mathrm T,1})_u$ due to Poincar\'e duality. Hence we have
$$\lim_{\mathcal P_{n_u} \ni \nu \to \infty} \frac{\mathrm{Tr}(F^{\nu}\circ \wt w,\mathscr H^{2d_u}(\mathscr R_{\mathrm T,1})_u)}{q^{\nu d_u}}=\frac{1}{q^{d_u}} \cdot \mathrm{Tr}(F\circ \wt w,\mathscr H^{2d_u}(\mathscr R_{\mathrm T,1})_u).
$$
 For $0\leq i \leq 2d_u-1$, we have
$$
\lim_{\mathcal P_{n_u} \ni \nu \to \infty} \frac{\mathrm{Tr}(F^{\nu}\circ \wt w,\mathscr H^{i}(\mathscr R_{\mathrm T,1})_u)}{q^{\nu d_u}}=0,
$$
by Theorem \ref{dls}.
It is clear that we have $\mathscr H^i(\mathscr R_{\mathrm T,1})_u= 0$ for $i \geq 1+2d_u$ or $i \leq-1$. 
Note that for $\nu \in \mathcal P_{n_u}$, the class of $T$ as a $F^\nu$-stable torus is $[N_{\mathrm W_G}^\nu (w)]$, where
$$N_{\mathrm W_G}^\nu (w)= w \cdot F(w) \cdot F^2(w) \cdot \ldots \cdot F^{\nu-1}(w). 
$$
We see that $N_{\mathrm W_G}^\nu(w)=w$ by the construction of $n_u$.
 Then we use Theorem \ref{dls} and Proposition \ref{ts} to represent
 $$Q_{T}^{G,\nu}(u)=R_{T,1}^{G,\nu}(u)=\mathrm{Tr}(F^\nu, (\mathscr R_{T,1})_u)=
\mathrm{Tr}(F^\nu \circ\wt{N_{\mathrm W_G}^\nu(w)},(\mathscr R_{\mathrm T,1})_u)=\mathrm{Tr}(F^\nu\circ \wt w,(\mathscr R_{\mathrm T,1})_u).$$
This completes the proof.
\end{proof}

\end{document}